\documentclass[12pt]{amsart}
\usepackage{amsthm}  
\usepackage{amssymb,mathrsfs,comment}
\usepackage{color,mathtools, mathabx}
\usepackage{tikz,tikz-cd,hyperref}
\usepackage{setspace}
\usepackage{tcolorbox} 

\usepackage{amsrefs} 

\usetikzlibrary{matrix,arrows,decorations.pathmorphing}

\usepackage{enumerate}

\theoremstyle{plain}
\newtheorem{theorem}{Theorem}[section]
\newtheorem{corollary}[theorem]{Corollary}
\newtheorem{lemma}[theorem]{Lemma}
\newtheorem{proposition}[theorem]{Proposition}

\theoremstyle{definition} 
\newtheorem{definition}[theorem]{Definition}

\newtheorem{remark}[theorem]{Remark}
\newtheorem{question}[theorem]{Question}
\numberwithin{equation}{section}

\title{Coloured Isomorphism of Classifiable C$^*$-algebras}
\date{}
\author{George A. Elliott and Jeffrey Im}

\begin{document}

\maketitle

\begin{abstract} 
  It is shown that the coloured isomorphism class of
  a
  unital, simple, $\mathcal{Z}$-stable, separable amenable 
  C$^*$-algebra
  satisfying the Universal Coefficient Theorem (UCT) is
  determined by its tracial simplex.
\end{abstract}

\section*{Introduction}

Approximate intertwinings have played a significant role in the
classification theory of C*-algebras (e.g., \cite{bratteli},
\cite{elliott_af}, \cite{pimsner_voiculescu},
\cite{elliott_loring}, \cite{elliott_evans}, \cite{elliott},
\cite{lin_tracial}, \cite{phillips}, \cite{kp}, 
\cite{de},
\cite{gong}, \cite{lin_tori},  
\cite{egl}, \cite{gln_1}, \cite{gln}, \cite{egln_kk}, \cite{gl}).
A related notion called $n$-coloured isomorphism, which
specializes to an approximate
intertwining in the case of a single colour,
was considered by Jorge Castillejos in
\cite{castillejos}. 
There it is shown that any two (unital) classifiable
C$^*$-algebras (i.e., those satisfying the hypotheses of Theorem
\ref{thm}) with at most one trace are 2-coloured
isomorphic, and it is posed as a question whether any two
such C$^*$-algebras
with isomorphic tracial simplices are
$n$-coloured isomorphic. We show that this question has an
affirmative answer using a somewhat modified notion, which we
shall
simply refer to as coloured isomorphism.

Let us begin by giving context for the notion from
\cite{castillejos} and comparing the approaches taken there and
here.   In fact, much of the original strategy is
  retained in this paper, so let us outline the
  proof of the main result, the
  finite case considered in Theorem \ref{thm},
  and mention the differences with \cite{castillejos} along the way. 
Throughout,
let $\omega$ be a fixed free ultrafilter on the natural
numbers. 
Two unital C$^*$-algebras $A$ and $B$ were said in
\cite{castillejos}
to be {\it $n$-coloured isomorphic} if there exist c.p.c.~ order zero
maps $\varphi : A \to B$ and $\psi : B \to
A$, and unitaries $u_1,
\ldots, u_n \in A_\omega$ and $v_1, \ldots, v_n \in B_\omega$, such
that
\[
  \sum_{k=1}^n u_k \psi \varphi(a) u_k^* = a
  \hspace{15 pt} 
  \text{ and }
  \hspace{15 pt} 
  \sum_{k=1}^n v_k \varphi \psi(b) v_k^* = b
  \]
  for all $a \in A$ and all $b \in B$.
  The present notion, coloured isomorphism, on the other hand, has
  order zero maps at the level of ultrapowers $A_\omega$ and
  $B_\omega$ (and so includes the case that the ultrapower maps
  are induced by a sequence of order zero maps at the level of the
  algebras $A$ and $B$); and the
  $u_i$ (resp.~ $v_i$) are contractions (not just unitaries)
  such that the absolute values squared
  add up to the identity of $A_\omega$ (resp.~
  $B_\omega$) (rather than a multiple of this). 
  (We also assume that traces on the ultrapower determined by a
  single trace on the algebra are preserved.)

  While the earlier notion, $n$-coloured isomorphism,
  preserves the tracial
  cone, up to isomorphism, the present notion preserves
  the tracial simplex, up to isomorphism (Theorem
  \ref{preserve_traces}).
  We note that isomorphism of 
  tracial cones coincides with (a multiple of) tracial simplex
  isomorphism in
  the cases considered in \cite{castillejos} (at most one trace).
  Our definition is formulated in terms of
  coloured equivalence for order zero maps,
  which has its origins in the
  $\mathcal Z$-stable implies finite nuclear dimension direction
  of the Toms-Winter conjecture, having first appeared
  in \cite{ms} and later more definitively in \cite{sww} and
  \cite{bbstww}. (That of \cite{castillejos} is not explicitly
  based on coloured equivalence of maps.)

  Now let $A$ and $B$ be finite classifiable C$^*$-algebras with
  isomorphic tracial simplices. In order to show that $A$ and $B$
  are coloured isomorphic, we must construct order zero maps
  $\varphi : A_\omega \to B_\omega$ and $\psi : B_\omega \to
  A_\omega$ such that $\psi \varphi \iota_A$ is coloured
  equivalent to $\iota_A$ where $\iota_A$ is the canonical
  embedding of $A$ into $A_\omega$, and likewise for $\varphi \psi
  \iota_B$ and $\iota_B$. By the finite order zero uniqueness
  theorem of \cite{bbstww}, which was later extended to remove 
  restrictions on the tracial simplices in \cite{cetww},
  it is sufficient to show that (in the notation of Corollary
  \ref{functional_calculus}) 
  $\tau (\psi \varphi)^n = \tau$ for all $n \in
  \mathbb N$ and all 
  $\tau \in T(A_\omega)$, and that $\tau(\varphi \psi)^n =
  \tau$
  for all $n \in \mathbb N$ and all $\tau \in T(B_\omega)$.
  Since $(\psi \varphi)^n = \psi^n \varphi^n$ and $(\varphi
  \psi)^n = \varphi^n \psi^n$, by Corollary
  \ref{multiplicative}, these tracial identities imply that
  $\varphi^n$ and $\psi^n$ also (for each $n$) induce mutually inverse isomorphisms of
  $T(A_\omega)$ and $T(B_\omega)$. In
  fact, the order zero maps we construct will induce the same
  mutually inverse isomorphisms of tracial simplices for each $n$.

  Roughly speaking, $\varphi$ is constructed with a sequence of
  maps at the level of the algebras $A$ and $B$ which factor
  through a fixed AF algebra $D$ via a fixed embedding
  $\alpha_A$. Furthermore, $\alpha_A$ is chosen so that it induces
  an isomorphism of the tracial simplices $T(D)$ and $T(A)$.
  This is obtained from the recently established
  homomorphism theorem of \cite{gln}. Such embeddings are now also
  known to exist outside the classification setting (i.e.,
  without the assumption of $\mathcal Z$-stability), building on
  ideas from \cite{sch17} and \cite{sch20}.
  Lastly, the maps $\varphi_k$
  into $B$ from the AF algebra $D$ are chosen to have prescribed
  tracial data. More specifically, given a faithful trace $\mu_k$
  on $C_0(0,1]$ and an affine map $\Phi : T(B) \to T(D)$,
  $\varphi_k$ is an order zero map satisfying the identity
  \[
    \tau \varphi_k^n = \mu_k(t^n) \Phi(\tau)
    \tag{$*$}
  \]
  for each $k, n \in \mathbb N$ and each $\tau \in T(B)$.
  (But, for obvious reasons, with $\mu_k(t^n) \to 1$ for each $n$.)
  This is the content of Theorem
  \ref{induce_tensor_trace}.

  Let us show how
  the desired order zero maps
  $\varphi$ and $\psi$ can be constructed from
  here. Let $\Phi$ be an isomorphism of the tracial
  simplices $T(B)$ and $T(D)$, which exists since $T(A)$ is
  assumed to be isomorphic to $T(B)$, by hypothesis, and $D$
  was chosen so that $T(D)$ is isomorphic to $T(A)$; and let
  $\mu_k$ be a sequence of faithful traces on $C_0(0,1]$ such that
  $\lim_{k \to \omega} \mu_k(t^n) = 1$ for each $n \in \mathbb
  N$, where $t$ denotes the identity map on $(0,1]$.
  (It is enough that $\lim_{k \to \infty} \mu_k(t) = 1$.)
  Then with $\varphi_k : D \to B$
  satisfying $(*)$, and with $\varphi: A_\omega \to B_\omega$
   the order zero map induced by the sequence 
  $(\varphi_k \alpha_A)$, we have
  \[
    \tau \varphi^n
    = \tau((\varphi_k \alpha_A)^n) =
     \tau(\varphi_k^n \alpha_A)
  \]
  for each $n \in \mathbb N$ and for each trace $\tau \in
  T(B_\omega)$ 
  with Corollary \ref{multiplicative} having been
  used at the 
  last equality.
  In order to make use of the fact that 
  $\mu_k$ levels out moments, we will need a reduction which was
  observed in \cite{ozawa}: The limit traces in $T(B_\omega)$
  (i.e., those which are induced by a sequence of traces in $T(B)$)
  are weak* dense in
  $T(B_\omega)$ (Theorem \ref{limit_traces}, below).
  Therefore, it suffices to check that the desired
  tracial identities (see next paragraph) hold
  for such traces.
  Let a trace $\tau$ of the form
  $\lim_{k \to \omega} \tau_k$ in $T(B_\omega)$ be given. Then
  \begin{align*}
    \tau(\varphi_k^n \alpha_A)
    & = \lim_{k \to \omega} \tau_k \varphi_k^n \alpha_A
      = \lim_{k \to \omega} \alpha_A^* \tau_k \varphi_k^n 
      = \alpha_A^* \lim_{k \to \omega} \tau_k \varphi_k^n \\
    &\stackrel{(*)}{=} \alpha_A^* \lim_{k \to \omega}
      \mu_k(t^n) \Phi(\tau_k)
      = \alpha_A^* \lim_{k \to \omega} \Phi(\tau_k)
    = \lim_{k \to \omega} \alpha_A^* \Phi(\tau_k)
  \end{align*}
  for each $k, n \in \mathbb N$.
  The second to last equality uses
  the fact that $\lim_{k \to \omega} \mu_k(t^n) = 1$ and
  continuity of $\alpha_A^*$ is used for the third equality and
  the last two.
  
  The remaining map $\psi$ is
  constructed in a similar way. Let $\alpha_B : B \to E$ be a
  AF embedding which induces a tracial simplex isomorphism $T(E)
  \to T(B)$, and
  let $\psi_k : E \to A$ be order zero maps satisfying the
  identity
  \[
    \tau \psi_k^n = \mu_k(t^n) \Psi(\tau)
  \]
  for each $k, n \in \mathbb N$ and each $\tau \in T(A)$, where
  $\Psi$ is the tracial simplex isomorphism
  $(\alpha_A^* \Phi \alpha_B^*) ^{-1}$. Let
  $\psi : B_\omega \to A_\omega$
  denote the order zero map induced by the sequence
  $(\psi_k \alpha_B)$.
  Then for each limit trace $\tau =
  \lim_{k \to \omega} \tau_k$ in $T(A_\omega)$,
  \[
    \tau \psi^n = \lim_{k \to \omega} \alpha_B^* \Psi(\tau_k)
  \]
  for each $n \in \mathbb N$. Therefore (see proof of Theorem
  \ref{thm} for more details),
  \[
    \tau (\psi \varphi)^n = 
    (\varphi^n)^*(\lim_{k \to \omega} \alpha_B^* \Psi(\tau_k))
    = \lim_{k \to \omega} \alpha_A^* \Phi \alpha_B^* \Psi(\tau_k)
    = \tau
  \]
  for each $n \in \mathbb N$ and each limit trace $\tau \in
  T(A_\omega)$. By Theorem \ref{limit_traces},
  the above identity holds for all $\tau \in T(A_\omega)$.
  A symmetric argument shows that
  $\tau (\varphi \psi)^n = \tau$ for each $n \in \mathbb N$ and
  each $\tau \in T(B_\omega)$. Therefore, $\varphi$ and $\psi$
  determine a coloured isomorphism of $A$ and $B$. (To simplify
  the discussion, we omit the question of preserving constant
  limit traces and defer to the proof of Theorem \ref{thm}.)

  Since the $n$-coloured isomorphism of \cite{castillejos} 
  requires unitaries rather than
  contractions, a different uniqueness theorem is developed
  for order zero maps in the unique trace case (\cite[Lemma
  5.6.2]{castillejos}). It provides a stronger 
  statement than \cite[Theorem 5.5]{bbstww} (in
  the unique
  trace case) since it is applicable to pairs of order zero maps
  (rather than one $*$-homomorphism and one order zero map) and
  because it provides unitary equivalence of the order zero maps
  involved rather than after tensoring the order zero maps with a
  positive element $h \in \mathcal Z$ with full spectrum.
  Let $h$ be such an element with the additional stipulation that
  $\tau_{\mathcal Z}(h^n) = \tau_{\mathcal Z}( (1_{\mathcal
    Z}-h)^n) = 1/(n+1)$ for each $n \in \mathbb N$
  where $\tau_{\mathcal Z}$
  denotes the unique trace of $\mathcal Z$.
  Because the scaling factors introduced by the unitaries
  differ from those introduced by contractions
  under traces, the compositions of the order zero maps $\varphi :
  A \to B$ and $\psi : B \to A$ implementing the $n$-coloured
  equivalence are compared with the contractive order zero maps
  $\rho_{A,h} := \sigma_1(\text{id}_A \otimes h)$ and 
  $\rho_{A,h} := \sigma_2(\text{id}_B \otimes h$) under traces 
  where $\sigma_1 : A \otimes \mathcal
  Z \to A$ and $\sigma_2 : B \otimes \mathcal Z \to B$ 
  are isomorphisms whose inverses are approximately unitarily
  equivalent to the first factor embeddings. 
  The moment problem in this setting is
  more complicated than the one
  in ours
  because the order zero maps involved
  in the uniqueness theorem are induced by constant sequences of
  maps, and therefore
  the moments need to match up on the dot rather than only
  approximately. Moreover, the target moments of $\psi
  \varphi$ and $\varphi \psi$ actually depend on $n \in \mathbb
  N$, as 
  \[
    \tau_A \rho_1^n  = \frac{\tau_A(\cdot)}{n+1}
    \hspace{10 pt}
    \text{ and } \hspace{10 pt}
    \tau_B \rho_2^n  = \frac{\tau_B(\cdot)}{n+1},
  \]
  where $\tau_A$ and $\tau_B$ denote the unique trace of $A$ and
  $B$ (see \cite[Theorem 5.6.8]{castillejos}).

  The main step in constructing $\varphi$ and $\psi$ (in
  \cite{castillejos}) is the
  construction of maps out of unital, simple, separable AF
  algebras with unique trace into $\mathcal Z$ realizing specific
  moments. Let us outline the construction of $\varphi$ 
  using this result. 
  By \cite[Theorem A]{sch17}, $A$ embeds into a
  simple, separable, unital AF algebra $D$ with unique trace via
  a map $\alpha_A$ which induces an isomorphism of tracial
  simplices.
  Then an order zero map 
  $\phi : D \to \mathcal Z$ is constructed, using a
  measure $\mu$ on the unit interval with the very specific
  moments
  $\mu(t^n) = 1/\sqrt{n+1}$,
  $n \in \mathbb N$,
  such that
  \[
    \tau_{\mathcal Z} \phi^n = \frac{\tau_D}{\sqrt{n+1}}
  \]
  for each $n \in \mathbb N$ (\cite[Lemma 5.6.5]{castillejos}). 
  The order zero map $\varphi$ is then
  given by the composition 
  \[
    \begin{tikzcd}
      A \arrow[r, "\textstyle \alpha_A", hook] & D \arrow[r,
      "\textstyle \phi"] & \mathcal{Z} \arrow[rr, "\textstyle 1_B
      \otimes \text{id}_{\mathcal Z}"] & & B \otimes
      \mathcal{Z} \arrow[r, "\textstyle \sigma_2"] & B.
    \end{tikzcd}
  \]
  Apart from $\phi$, each of the above maps is a trace-preserving
  $*$-homomorphism, and so by Corollary \ref{multiplicative}
  below, the
  moments of $\varphi$ are
  \[
    \tau_B \varphi^n = \tau \sigma_2 (1_B \otimes 1_{\mathcal Z})
    \phi^n \alpha_A = \frac{\tau_B}{\sqrt{n+1}}
  \]
  for each $n \in \mathbb N$. The second order zero map $\psi$
  is constructed in a similar way so that 
  \[
    \tau_A (\psi \varphi)^n = \frac{\tau_A}{n+1}
    \hspace{10 pt} \text{ and }
    \hspace{10 pt} 
\tau_B (\varphi \psi)^n = \frac{\tau_B}{n+1}
  \]
  for each $n \in \mathbb N$.
  The appropriate
  uniqueness theorem (\cite[Lemma 5.6.2]{castillejos}) is then applied twice
  with $\psi \varphi$: once with $\rho_{A,h}$, and once more with
  $\rho_{A,1-h}$, and similarly for $\varphi \psi$, in order to
  obtain a two-coloured isomorphism of $A$ and $B$.

  All that remains to be outlined in the construction of order
  zero maps, either in our setting or in that of
  \cite{castillejos}, is how to ensure the
  prescribed tracial data.
  Let us begin with the approach taken in
  \cite{castillejos}.
  Let $D$ be a simple, separable, unital AF algebra with unique
  trace $\tau$ and let $\nu$ be a fully supported Borel measure on
  $[0,1]$. Recall (e.g., from \cite{ers}) that $\tau$ induces a
  functional $d_\tau : \text{Cu}(D) \to [0,\infty]$. 
  It is stated in \cite[Proposition 1.10.12]{castillejos}
  that the map
  $\sigma: \text{Lsc}([0,1],\text{Cu} (D)) \to \text{Cu}
  (\mathcal Z)$ 
  determined by the rule
  \[
    f \mapsto \int_0^1 d_\tau (f(t)) ~ d\nu(t)
  \]
  for each $f \in \text{Lsc}([0,1],\text{Cu} (D))$ is a Cuntz
  category morphism. Since $\mathcal Z$ has unique trace,
  $\sigma$ maps into $\text{Cu} (\mathcal Z)$ (which is naturally
  isomorphic to $V(D) \sqcup (0,\infty]$, by \cite[Corollary
  6.8]{ers}). By \cite[Theorem 2.6]{aps}, the natural map
  $\text{Cu}(C([0,1],D)) \to \text{Lsc}([0,1],\text{Cu} (D))$ is
  a Cuntz category isomorphism. Combining this isomorphism
  with $\sigma$,
  one has a Cuntz category morphism from $\text{Cu}(C([0,1],D))
  \to \text{Cu}(\mathcal Z)$. This induces a Cuntz category
  morphism between $C_0(0,1] \otimes D$ and $\mathcal Z$  
  of the augmented invariant of \cite{robert} and by
   \cite[Theorem 1.0.1]{robert}, there exists a
  $*$-homomorphism $\pi : C_0(0,1] \otimes D \to \mathcal Z$ which
  then induces an order zero map $\rho : D \to \mathcal Z$, by
  \cite[Corollary 4.1]{wz_cpc}. It is then checked 
  that $\rho$ satisfies the desired tracial identity in
  \cite[Lemma 5.6.5]{castillejos}.
  
  In order to move beyond the unique trace case in Theorem
  \ref{thm}, the order zero maps implementing a coloured
  equivalence cannot be made with a sequence of maps that factor
  through $\mathcal Z$, as that would collapse the trace space to
  a single point (such maps necessarily preserve the trace space
  up to isomorphism, by Theorem \ref{preserve_traces} below).
  Hence, a replacement for \cite[Lemma
  5.6.5]{castillejos} is needed. 
  Our formulation of the order zero map realizing prescribed
  tracial data maps directly into a unital, simple, separable,
  exact, $\mathcal Z$-stable C$^*$-algebra $B$ with stable rank
  one whose tracial simplex is isomorphic to that of an AF algebra
  $D$ (no longer assumed to have unique trace) instead of
  $\mathcal Z$. 
  This order zero map is constructed by showing that 
  the map $\sigma : \text{Cu}(C_0(0,1] \otimes D)
  \to \text{LAff}_+(TB)$
  determined by the rule
  \[
    [d] \mapsto (\mu \otimes \Phi(\cdot))[d]
  \] for each $[d] \in \text{Cu}(C_0(0,1] \otimes D)$
  is a Cuntz category morphism, where $\mu$ is a faithful densely
  defined lower semicontinuous trace on $C_0(0,1]$.
  Since $\text{LAff}_+(TB)$ is a
  subobject of $\text{Cu}(B)$ (Remark \ref{subobject}),
  $\sigma$ extends to a Cuntz category morphism
  into $\text{Cu}(B)$.

  In fact, we couldn't figure out how to use
  the characterization of compact containment given in \cite{aps}
  to work directly at the level of the cone over the AF algebra
  and so it is established that $\sigma$ is a Cuntz category
  morphism by showing that it is the inductive limit of Cuntz
  category morphisms $\sigma_i$ at the finite stages of the AF
  algebra inductive limit decomposition.  That
  $\sigma_i$ is a generalized Cuntz morphism follows from
  \cite{ers}. To show that compact containment is preserved by
  $\sigma_i$, we show (in the proof of Theorem
  \ref{induce_tensor_trace}) that $\sigma_i$ is
  a weighted
  direct sum of copies of
  the functional $d_\mu : \text{Cu}(C_0(0,1]) \to
  [0,\infty]$. This reduces the problem of showing that $\sigma_i$
  preserves compact containment to showing that the functional
  $d_\mu$ preserves compact containment. This is done in Section
  \ref{functionals}, where we give a topological
  characterization of when functionals arising from a faithful
  densely defined lower semicontinuous trace preserve compact
  containment. In particular, an essential property of the half open
  interval is that it does not contain non-empty compact open sets.
  (At the end of the section, we give
  a sufficient condition for when non-faithful densely defined
  lower semicontinuous traces induce functionals which preserve
  compact containment.) In
  Theorem \ref{induce_tensor_trace}, it is shown that
  the morphisms $\sigma_i$ give rise to a one-sided intertwining
  at the level of the Cuntz category. By the 
  classification theorem (of \cite{ce})
  for $*$-homomorphisms out of cones, this then
  gives rise to an approximate one-sided intertwining at the level
  of C$^*$-algebras. The inductive limit $*$-homomorphism $\pi :
  C_0(0,1] \otimes D \to B$ then induces the desired order zero
  map.

  The outline of the paper is as follows: In Section \ref{preliminaries}, we
  introduce basic notions and terminology along with important
  results that are essential for the main results.
  We also include statements which are likely known to
  experts, but possibly not spelled out in
  literature. In Section \ref{functionals}, we give a
  characterization for when a trace $\tau$ on a
  C$^*$-algebra $A$ induces a functional $d_\tau$ on $\text{Cu}(A)$
  preserves compact containment.
  In Section \ref{main_section}, we prove the main results and
  discuss some consequences and questions.
  
  \subsection*{Acknowledgements}
  
  The second author would like to thank Jorge Castillejos,
  Christopher Schafhauser, Aaron Tikuisis, and Stuart White for
  suggestions and encouragement.

\section{Preliminaries} \label{preliminaries}

\subsection{Order zero maps}
Let $A$ and $B$ be C$^*$-algebras and let $\varphi : A \to B$ be
a completely positive (c.p.) map (c.p.c.~if the map is
contractive). The map $\varphi$ is said to have {\it
  order zero} if it preserves orthogonality ({\cite[Definition
  1.3]{wz_cpc}}).
Examples of such
maps include $*$-homomorphisms and more generally,
the product of a $*$-homomorphism $\pi : A \to B$
with a positive element $h \in B$
which commutes with the image of $\pi$. The following structure
theorem shows that all order zero maps admit such a
decomposition. The $*$-homomorphism $\pi_\varphi$ which appears
below is called the support $*$-homomorphism of $\varphi$.
\begin{theorem}[{\cite[Theorem 2.3]{wz_cpc}}] \label{structure}
  Let $A$ and $B$ be
  C$^*$-algebras and let $\varphi : A \to B$ be a c.p.~ order zero
  map. Define $C := \text{C}^*(\varphi(A)) \subseteq B$
  (and denote the multiplier algebra of $C$ by $\mathcal M(C)$). 
  Then there
  exist a unique positive element $h_\varphi$ in the centre of
  $\mathcal M(C)$
  and a $*$-homomorphism $\pi_\varphi : A \to
  \mathcal M(C)$ such that
  \begin{align}
    \varphi(a) = h_\varphi  \pi_\varphi(a)
    \label{structure_eq}
  \end{align}
  for all $a \in A$. Necessarily, $\|h_\varphi\| = \|\varphi\|$, and
  if $A$ is unital, then $h_\varphi = \varphi(1_A) \in
  B$.
\end{theorem}

\begin{corollary}[{\cite[Proposition 1.4]{bbstww}}] \label{projection}
  Let $\varphi : A \to B$ be a c.p.c.~ order zero map between
  C$^*$-algebras where $A$ is unital. Then $\varphi$ is a
  $*$-homomorphism if, and only if, $\varphi(1_A)$ is a projection.
\end{corollary}

\begin{proof}
  This follows immediately from Theorem \ref{structure}.
\end{proof}

The next corollary is known as the order
zero functional calculus.
To avoid ambiguity with the notation
introduced below, we will never use $\varphi^n$ to mean
iterated composition. 

\begin{corollary}[{\cite[Corollary 3.2]{wz_cpc}}] Let $\varphi :
  A \to B$ be a c.p.c.~ order zero map and let notation be as in
  Theorem \ref{structure}. Then for any positive function $f$
  in $C_0((0,1])$, the map
  \[
    f(\varphi) : A \to B
  \]
  defined by
  \begin{align}
    f(\varphi)(\cdot) := f(h_\varphi)\pi_\varphi(\cdot)
    \label{fc}
  \end{align}
  is a well-defined c.p.~ order zero map taking values in $C$. If
  $f$ has norm at most one, then the map
  $f(\varphi)$ is contractive.
  \label{functional_calculus} 
\end{corollary}

Let $\varphi : A \to B$ and $\psi : B \to C$ be c.p.c.~ order zero
maps between unital C$^*$-algebras and let notation be as in
Theorem \ref{structure}.
The following statement concerning the composition $\psi \varphi$
is  
established in \cite{castillejos} by proving that
the support $*$-homomorphism of a composition is essentially the
composition of the the individual support $*$-homomorphisms (more
precisely, it is shown that
$h_{\psi
  \varphi}\pi_{\psi \varphi} = h_{\psi \varphi} \pi_\psi \pi_\varphi$
(\cite[Corollary 1.4.14]{castillejos})). We give a slightly more
direct proof.
\begin{corollary}[{\cite[Section 5.6]{castillejos}}] \label{multiplicative}
  Suppose $\varphi : A \to B$ and $\psi : B
  \to C$ are c.p.c.~ order zero maps between C$^*$-algebras.
  Then
  \begin{align}
    (\psi \varphi)^n = \psi^n \varphi^n
    \label{multiplicative_id}
  \end{align}
  for each $n \in \mathbb N$ (excluding $n=0$). If $\varphi$ is a
  $*$-homomorphism, then $\varphi^n = \varphi$.
\end{corollary}

\begin{proof}
  By \cite[Proposition 2.2]{wz_cpc}, 
  we may suppose, without loss of generality, that $A$ is
  unital.
  Let notation be as in Theorem \ref{structure} for $\varphi$,
  $\psi$, and $\psi \varphi$. Then we have
  \begin{align}
    h_{\psi \varphi} = \psi \varphi(1_A) = \psi(h_\varphi)
    \stackrel{(\ref{structure_eq})}{=} h_\psi \pi_\psi(h_\varphi),
    \label{positive_multiplicative}
  \end{align}
  using that $A$ is unital for the first two equalities. Now,
  for fixed $n \in \mathbb N$, we have
  \begin{align*}
    (\psi \varphi)^{n+1}
    & \stackrel{(\ref{fc})}{=} h_{\psi \varphi}^{n+1} \pi_{\psi
      \varphi}
      = h_{\psi \varphi}^n h_{\psi \varphi} \pi_{\psi \varphi}      
     \stackrel{(\ref{positive_multiplicative})}{=} 
      h_\psi^n \pi_\psi(h_\varphi^n) h_{\psi \varphi} \pi_{\psi
      \varphi}
    \\
    &      \stackrel{(\ref{structure_eq})}{=} h_{\psi}^n
  \pi_\psi(h_\varphi^n) \psi \varphi
      \stackrel{(\ref{structure_eq})}{=} (h_{\psi}^{n+1}
      \pi_\psi(h_\varphi^n)\pi_\psi)(h_\varphi \pi_\varphi) 
    \\ & 
      = (h_\psi^{n+1} \pi_\psi)(h_\varphi^{n+1} \pi_\varphi)
      \stackrel{(\ref{fc})}{=} \psi^{n+1} \varphi^{n+1}.
  \end{align*}
  If $\varphi$ is a $*$-homomorphism, then $h_\varphi = \varphi(1_A)$ is
  a projection and so $\varphi^n = \varphi$.
\end{proof}

The following correspondence between c.p.c.~ order zero maps and
$*$-homomorphisms out of cones
is from \cite{wz_cpc}, but the formulation
presented here is that of \cite[Proposition 1.3]{bbstww}. 
Both
Corollary \ref{cone_oz} and Remark \ref{cone_fc} will be used
in the proof of Theorem \ref{thm}.

\begin{corollary}[{\cite[Corollary 3.1]{wz_cpc}}]
  \label{cone_oz}
  Let $A$ and
  $B$ be C$^*$-algebras. There is a one-to-one correspondence
  between c.p.c.~ order zero maps $\varphi : A \to B$ and
  $*$-homomorphisms $\pi: C_0(0,1] \otimes A \to B$ where
  $\varphi$ and $\pi$ are related by the commutating diagram
  \begin{equation}
    \begin{tikzcd}
      A \arrow[rr, "\textstyle a \mapsto t \otimes a"] \arrow[rrd,
      "\textstyle \varphi"'] &  & {C_0(0,1] \otimes A} \arrow[d,
      "\textstyle \pi"] \\
      &  & B
    \end{tikzcd}
    \label{cone_diagram}
  \end{equation}
  and $t \in C_0(0,1]$ denotes the identity function.
\end{corollary}

\begin{remark} \label{cone_fc}
  It is noted in \cite{bbstww} that the order zero functional
calculus can be recovered from the identity
\begin{align} 
  f(\varphi)(a) = \pi(f \otimes a)
  \label{cone_eq}
\end{align}
where
$f$ is a positive function in $C_0(0,1]$, $a \in A$,
and $\varphi$ and $\pi$
are as in the corollary above. Let us verify that this agrees with
the definition made earlier.

\begin{proof}
  We may again, by \cite[Proposition 2.2]{wz_cpc}, suppose without
  loss of generality that $A$ is unital. Let notation be as in
  Theorem \ref{structure} and let $f \in C_0(0,1]_+$ be given.
  By Theorem \ref{structure} and Corollary \ref{cone_oz},
  $h_\varphi = \varphi(1_A) = \pi(t \otimes 1_A)$. So for fixed $n
  \in \mathbb N$ and $a \in A$,
  \begin{align*} 
    \varphi^{n+1}(a) & \stackrel{(\ref{fc})}{=} h_\varphi^{n+1} \pi_\varphi(a) =
                    h_\varphi^n (h_\varphi \pi_\varphi(a)) \stackrel{(\ref{structure_eq})}{=}
                    h_\varphi^n \varphi(a) \\
                  & \stackrel{(\ref{cone_diagram})}{=} \pi(t^n \otimes 1_A) \pi(t \otimes a)
                    =  \pi(t^{n+1} \otimes a).
  \end{align*}
  It is readily seen from this calculation and
  approximating $f$ by polynomials in $C_0(0,1]$
  that $f(\varphi)(a) = \pi(f \otimes a)$.
\end{proof}
\end{remark}

\subsection{Ultrapowers} Let $A$ be a
C$^*$-algebra. The bounded sequence algebra of $A$, denoted by
$l^\infty(A)$, is defined to be the collection of all
norm-bounded sequences of elements from $A$. The ideal of elements in
$l^\infty(A)$ which go to zero along the ultrafilter $\omega$ is
denoted by $c_\omega(A)$. The {\it ultrapower} of $A$ is the
quotient
C$^*$-algebra $A_\omega := l^\infty(A)/c_\omega(A)$. We shall write
$(a_n)_{n=1}^\infty$ for an element in $A_\omega$ rather than the
class it belongs to. We will occasionally
write $a$ to refer to the image of $a$ under
the canonical embedding of $A$ into $A_\omega$. A sequence $\varphi_n : A
\to B$ of c.p.c.~ order zero maps between
C$^*$-algebras induces a c.p.c.~ order zero map $\varphi :
A_\omega \to B_\omega$ between ultrapowers. We denote the induced
map by $\varphi := (\varphi_n)_{n = 1}^\infty$. If $A$ is unital,
then, with notation as in Theorem \ref{structure},
we see that $h_\varphi = (h_{\varphi_n})_{n=1}^\infty$ and so the
order zero functional calculus for $\varphi$ can be realized by
applying positive functions in $C_0(0,1]$
componentwise to $\varphi_n$.
\subsection{Traces and finiteness} \label{traces}

Let $A$ be a C$^*$-algebra.
By a {\it trace} on $A$ we will mean a lower semicontinuous function $\tau
: A_+ \to [0,\infty]$ that is additive, preserves zero, is
positively
homogeneous, and satisfies the trace identity $\tau(a^*a) =
\tau(aa^*)$ for all $a \in A$.
We will
denote by $N_\tau$ the ideal of elements $a \in A$ such that
$\tau(a^*a) = 0$.
A trace $\tau$ is said to be {\it faithful} if 
$\tau(a^*a) = 0$ implies $a = 0$.
We will use the notation $T(A)$ (or simply $TA$) to
denote the collection of all tracial states on $A$ -- a Choquet
simplex if $A$ is unital. 
The {\it limit traces} of 
$A_\omega$ are the
tracial states on $A_\omega$ which are equal to
$\lim_{n \to \omega} \tau_n$ where
$(\tau_n)$ is some sequence of tracial states on $A$. We will denote the
collection of limit traces on $A_\omega$ by $T_\omega(A_\omega)$.
We will occasionally not make a notational distinction between
a trace in $T(A)$ and (what we will refer to as)
the constant limit trace in $T(A_\omega)$
induced by it.

By \cite[Corollary 3.4]{wz_cpc}, a c.p.c.~
order zero map $\varphi : A \to B$ induces a mapping of bounded
traces $\varphi^* : \mathbb R_+ T(B) \to \mathbb R_+ T(A)$.
We would be be interested in the case that
$\varphi^*$ is an 
isomorphism of tracial simplices, but this hope turns out not to
be realistic.
Rather, considering a sequence of c.p.c.~
order zero maps $\varphi_k : A \to B$, we will be interested in
the case that the c.p.c.~ order zero map
$\varphi : A_\omega \to
B_\omega$ induced by $(\varphi_k)$ preserves constant limit traces
(i.e., takes constant limit traces on $B_\omega$ to constant limit
traces on $A_\omega$)
and furthermore
the map $ \tau \mapsto \lim_{k \to
  \omega} \tau \varphi_k$, $\tau \in T(B)$, is an
isomorphism of the simplices $T(B)$ and $T(A)$.
(Note that requiring the sequence $(\varphi_k)$ to preserve
constant limit traces is equivalent 
to requiring the sequence $(\tau \varphi_k)$
to be norm convergent for every $\tau \in T(B)$.)
More generally, we shall be interested in the condition that
a constant limit trace preserving c.p.c.~ 
order zero map $\varphi : A_\omega \to B_\omega$ 
(not necessarily arising from a sequence $(\varphi_k)$ as above)
induces an isomorphism $T(B) \to T(A)$
of tracial simplices via the composed map
$\iota_A^* \varphi^* c_B$ 
where $c_B : T(B) \to T(B_\omega)$ denotes the
embedding of $T(B)$ as constant limit traces on
$T(B_\omega)$.  (One might also consider the tracial cones,
instead, and ask when they are isomorphic.)

Calculations involving traces on $A_\omega$
will often be reduced to the case of limit traces by
using the following fact about weak* density of limit traces.
Generalizations of the statement presented here can be found
in \cite{ozawa}, \cite{nr}, and \cite{rordam}.

\begin{theorem}[{\cite[Theorem 8]{ozawa}}]
  Let $A$ be a separable, exact, and $\mathcal
  Z$-stable C$^*$-algebra (where $\mathcal Z$ denotes the Jiang-Su
  algebra, \cite{jiang-su}). Then $T_\omega(A)$ is weak* dense in
  $T(A_\omega)$.
  \label{limit_traces}
\end{theorem}

Let $A$ be a C$^*$-algebra with nonempty tracial state space.
We define the seminorm $\| \cdot
\|_{2, T(A)}$ by
\[
  \|a\|_{2,T(A)} := \sup_{\tau \in T(A)} (\tau(a^*a))^{1/2},
\]
for each $a \in A$. The {\it trace-kernel ideal} of $A_\omega$ is
the set
\[
  J_A := \{(a_n) \in A_\omega :
  \lim_{n \to \omega} \|a_n\|_{2,T(A)} = 0
  \},
\]
and the quotient of $A_\omega$ by $J_A$ is called the {\it uniform
  tracial ultrapower} of $A$, and it is denoted by
$A^\omega$. More details about these notions
can be found in \cite[Section 4]{kirchberg_rordam}
and \cite[Section 1]{cetww}.

\begin{remark} \label{uniform_tracial_ultrapower}
  Let $A$ be a C$^*$-algebra and let $B$
  be a C$^*$-algebra with $T(B)$ nonempty and compact such that the limit traces of
  $B_\omega$ are dense in $T(B_\omega)$ (for example, unital
  C$^*$-algebras as in Theorem \ref{limit_traces})
  and let 
  $\varphi, \psi: A \to B_\omega$ be c.p.c.~
  order zero maps. Then $\tau \varphi = \tau \psi$ for all $\tau
  \in T(B_\omega)$ if, and only if, $\varphi$ agrees with $\psi$
  in the uniform tracial ultrapower $B^\omega$. 
\end{remark}

\begin{proof} Let $\pi : B_\omega \to B^\omega$ denote the
  canonical quotient map.
By linearity, it is enough to show that
an element of $B_\omega$ is in the kernel of every trace on
$T(B_\omega)$ exactly when it belongs to the trace-kernel ideal of
$B_\omega$. Suppose $(b_n) \in \ker \tau$ for every $\tau \in T(B_\omega)$.
Since $T(B)$ is compact, for each $n \in \mathbb N$, there exists
a trace $\tau_n \in T(B)$ such that $\|b_n\|_{2,T(B)}^2 =
\tau_n(b_n^*b_n)$. By assumption, $(b_n)$ is in the kernel
of the limit trace $\lim_{n \to \omega} \tau_n$ on
$T(B_\omega)$. Since the kernel of a trace is a left ideal, it
follows that
\[
  \|\pi(b_n)\|^2 = \lim_{n \to \omega}
\|b_n\|_{2,T(B)}^2 = \lim_{n \to \omega} \sup_{\tau \in
  T(B)}\tau(b_n^* b_n) = \lim_{n \to \omega} \tau_n(b_n^*b_n) =
0.
\]
Conversely, suppose $b = (b_n)$ is in the trace-kernel ideal of
$B_\omega$.
Since the limit traces of $B_\omega$ are weak* dense in $T(B_\omega)$, it is enough to show that
$\tau(b) = 0$ for every limit trace $\tau =
\lim_{n \to \omega} \tau_n \in T(B_\omega)$. This follows from the
computation
\[
  |\tau(b)| = \lim_{n \to \omega} |\tau_n(b_n)| \leq \lim_{n \to
    \omega} \tau_n(b_n^* b_n) \leq \lim_{n \to \omega}
  \|b\|_{2,T(B)}^2 = \|\pi(b)\| = 0.
\]
The first inequality follows from a well known fact about finite
traces (\cite[Theorem 1]{gardener}).
\end{proof}

The following statement holds somewhat more
generally (see proof), but this is the setting which
will be relevant for us in Theorem \ref{thm}. \vspace{10 pt} 

\begin{theorem}[\cite{rordam},\cite{cuntz_df},\cite{bh},\cite{haagerup},\cite{rieffel}]
  Let $A$ be a (non-zero)
  unital, simple, exact, $\mathcal{Z}$-stable
  C$^*$-algebra. Then $A$ is either purely infinite or finite.
  \label{dichotomy}
\end{theorem} 

\begin{proof}
  By \cite[Corollary 5.1]{rordam},
  the unital exact, $\mathcal Z$-stable C$^*$-algebra $A$ is
  purely infinite if, and only if, it is
  traceless (i.e. $T(A) = \varnothing$).
  By \cite[Corollary 4.7]{cuntz_df} and \cite[Theorem II.2.2]{bh},
  a simple, unital C$^*$-algebra $A$ is stably finite if,
  and only if, it admits a (non-zero, finite) quasitrace. 
By \cite[Theorem
  5.11]{haagerup}, every quasitrace on a unital
  exact C$^*$-algebra is a trace. 
  Therefore, $A$ is stably finite if, and only if, $T(A) \neq
  \varnothing$. 
  It remains to show that finiteness
  implies stable 
  finiteness in the present context.
  By \cite[Theorem 6.7]{rordam},
  since $A$ is simple,
  unital, $\mathcal Z$-stable, and finite, it has stable rank
  one.
  By
  \cite[Theorem 6.1]{rieffel}, $M_n(A)$ also has stable rank one,
  and, in particular, is finite
  (see \cite[Proposition 3.1]{rieffel});
  in other words, $A$ is stably finite.
\end{proof}

\begin{remark} \label{finiteness}
  When $A$ is unital, simple, and $\mathcal Z$-stable,
  the proof of the preceding theorem shows that the statement of
  \cite[Theorem 6.7]{rordam} can be slightly strengthened to
  say that the following three statements are equivalent:
  \begin{enumerate}[(1)]
  \item $A$ is finite;
  \item $A$ has stable rank one;
  \item $A$ is stably finite.
  \end{enumerate}
\end{remark}

\subsection{Coloured isomorphism} \label{colored}

We introduce two notions of coloured isomorphism. 
The first one (Definition \ref{coloured_isomorphic} below)
is symmetrically formulated, but it 
is a little long. A weaker notion (Definition \ref{minimalist}
below) is all that's needed to establish Theorem
\ref{preserve_traces}. Together with Theorem \ref{thm}, it follows
that these two notions
coincide in the classifiable setting.

\begin{definition} \label{coloured_isomorphic}
  Unital C$^*$-algebras $A$ and $B$ will be
  said to be {\it coloured isomorphic} 
  if there exist c.p.c.~order zero maps
  $\varphi : A_\omega \to B_\omega$ and $\psi : B_\omega \to
A_\omega$ such that $\psi \varphi \iota_A$ is coloured equivalent
to $\iota_A$ (see Definition \ref{coloured_equivalent}) and $\varphi \psi
\iota_B$ is coloured equivalent to 
$\iota_B$; and
$\varphi^*$ and $\psi^*$ preserve constant limit traces
(see Section \ref{traces}).
\end{definition}

\begin{remark}
  Coloured isomorphism is reflexive and
  symmetric. It follows from Corollary \ref{main_equivalence} that
  transitivity holds for
  classifiable C$^*$-algebras (cf.~ \cite[Section 6]{bbstww}),
  but transitivity is not clear in general.
\end{remark}

\begin{definition} \label{minimalist}
  Unital C$^*$-algebras $A$ and $B$ will be said to be
  {\it minimalist coloured isomorphic} if
  there exist constant limit trace preserving
c.p.c.~ order zero maps $\varphi: A_\omega \to
B_\omega$ and $\psi : B_\omega \to A_\omega$, and 
contractions $u_1, \ldots, u_m \in
A_\omega$ and $v_1, \ldots, v_n \in B_\omega$, such that
\begin{align}
  \sum_{i=1}^m u_i \psi \varphi(a) u_i^* = a
  \hspace{10 pt} \text{ and } \hspace{10 pt} 
  \sum_{j=1}^n v_j \varphi \psi (b) v_j^* = b
  \label{coloredid}
\end{align}
for all $a \in A$ and for all $b \in B$, and,
in addition,
\begin{align}
  \sum_{i=1}^m u_i^* u_i =
  1_{A_\omega}, \hspace{10 pt}
  \sum_{j=1}^n v_j^* v_j = 1_{B_\omega}.
  \label{coloredpartition}
\end{align}
\end{definition}
We note that (\ref{coloredid}) is not a priori symmetric,
as required by Definition \ref{coloured_isomorphic}.
It follows from Theorem \ref{preserve_traces} and Theorem
\ref{thm} that symmetry is automatic in the classifiable setting.

The identity (\ref{coloredid}) implies that
$\varphi$ and
$\psi$ are injective; furthermore,
the commutation relations coming from the coloured equivalences of
maps -- the symmetrized form of coloured isomorphism holding
automatically in the classifiable case -- guarantee that $u_i \psi
\varphi(\cdot) u_i^*$ and $v_j \varphi \psi(\cdot) v_j^*$ are
order zero $(i = 1, \ldots, m \text{ and } j = 1, \ldots,
n)$.

\subsection{The Cuntz category}

If $a$ and $b$ are elements of an ordered set $M$, we will say
that $a$ is {\it (countably) compactly contained} in $b$, and we write $a \ll
b$, if for any
increasing sequence $(b_k) _{k=1}^\infty$ in $M$ with $\sup_k b_k
\geq b$ (or such that every upper bound of the sequence majorizes $b$),
eventually $b_k \geq a$. An increasing sequence with each term compactly contained
  in the next is called {\it rapidly increasing}.
An element which is compactly contained in
itself is called {\it compact} (more precisely, {\it countably
  compact}).

  {\it Cuntz category
    semigroups} are ordered abelian
  semigroups with an additive identity with the following four
  properties:
  \begin{enumerate}
\item[(O1)] Every increasing sequence in $S$ has a supremum in $S$.
\item[(O2)] Any element of $S$ is the supremum of a rapidly increasing
  sequence.
\item[(O3)] If $a_i$ and $b_i$ are elements of $S$ such that
  $a_i \ll b_i$, $i = 1,2$, then $a_1 + a_2 \ll b_1 + b_2$.
\item[(O4)] If $(a_n)$ and $(b_n)$ are increasing sequences in $S$,
  then $\sup_n (a_n + b_n) = \sup_n a_n + \sup_n b_n$.
\end{enumerate}

{\it Cuntz category morphisms} $f : S \to T$ are ordered semigroup
maps (i.e. preserving order, addition, and the additive identity)
which preserve suprema of increasing sequences and compact
containment. A {\it generalized Cuntz category morphism} is a
Cuntz category morphism which does not necessarily preserve
compact containment.
The Cuntz category has as objects
  the ordered semigroups and ordered semigroup maps with the
  properties stipulated above (\cite{cei}, \cite{ers},
  \cite{robert}). It is easily checked, and it will be used
  without mention that Cartesian products and direct sums
  exist in this category.

  It was shown in \cite{cei} that there is a functor
  $\text{Cu}(\cdot)$
  from the category of C$^*$-algebras
  to the Cuntz category.
  Recall that, for a C$^*$-algebra $A$, $\text{Cu}(A)$ is
  the ordered
  semigroup of Cuntz equivalence classes of positive elements in
  the stabilization of $A$ and
  is a Cuntz category semigroup. We shall denote
  by $[a]$ the Cuntz equivalence class of a positive element $a$
  of $A \otimes \mathcal K$.
  We shall denote by $V(A)$ the semigroup of
Murray-von Neumann equivalence classes of projections in the
stabilization of $A$. For a compact convex subset $K$ of a locally
convex topological vector space, we shall denote by
$\text{LAff}_+(K)$ the collection of lower semicontinuous extended
positive real-valued affine functions on $K$ which are strictly
positive, except for the zero function, and are the pointwise
supremum of an increasing sequence of
continuous and finite-valued such functions. If $K$ is metrizable,
then the latter condition is automatic. It was observed in 
Section 2 of \cite{eik} (and in Section 3 of \cite{thiel}) that
$\text{LAff}_+(K)$, equipped with pointwise order and addition,
is a Cuntz category semigroup. The addition and order structure
defined on the disjoint union decomposition
appearing in the following important
computation can be found in Section 6 of \cite{ers}.
Related work on 
the unstabilized Cuntz semigroup that was done prior
to the result stated here can be found in \cite{pt},
\cite{bpt}, and \cite{dt}.
Recent developments that go well
beyond what is needed in this paper include \cite[Theorem
8.11]{thiel} and \cite[Theorem 7.15]{aprt}.

\begin{theorem}[{\cite[Corollary 6.8]{ers}}] \label{valaff}
  Let $A$ be a unital, simple, separable, exact, $\mathcal
  Z$-stable, and finite C$^*$-algebra. Then, in a natural way, 
  \[
    \text{Cu}(A) \cong
    V(A) \sqcup \text{LAff}_+(TA) \backslash \{0\}
    .
  \]
\end{theorem}

\begin{remark} \label{subobject}
  We note that $\text{LAff}_+(TA)$ is a subobject of
  $\text{Cu}(A)$
  in the setting of Theorem \ref{valaff}, according to the
  embedding there.

  \begin{proof}
    Since the elements of $V(A)$ are compact in
    $\text{Cu}(A)$ (by Corollary 5 of \cite{cei}; this statement
    does not use stable rank one), it
    is clear that the embedding of $\text{LAff}_+(TA)$ into
  $\text{Cu}(A)$ is a generalized Cuntz category morphism. To show
that compact containment is preserved, let $f \ll g$ in
$\text{LAff}_+(TA)$ and let $(g_n) _{n=1}^\infty$ be an increasing
sequence in $\text{Cu}(A)$ such that
$\sup_n g_n$ majorizes $g$. Passing to a subsequence, we may
suppose 
that the sequence $(g_n) _{n=1}^\infty$ is composed entirely of
projections or entirely of affine functions. In the latter case,
since the embedding of $\text{LAff}_+(TA)$ in $\text{Cu}(A)$
preserves increasing sequential suprema, we have $\sup_n g_n \in
\text{LAff}_+(TA)$ and so 
there is, by assumption, some $g_n$ which majorizes $f$. Let us now
consider the former case.
Denote by $\widehat g_k$ the rank of
$g_k$ -- the element of
$\text{LAff}_+(TA)$ got by evaluation of traces in $TA$ at the
projection $g_k$.

By the characterization of compact containment
in $\text{LAff}_+(TA)$
in Section 2 of \cite{eik} (or \cite[Lemma 3.6]{thiel}), there
exists an $h \in \text{LAff}_+(TA)$ which is continuous and
finite-valued and an $\varepsilon > 0$ such that $f \leq h < (1 +
\varepsilon) h \leq g$.
Since $0 \notin \text{LAff}_+(TA)$,
$h$ is strictly positive. 
It follows that there is, for each $\tau \in
TA$, a neighbourhood of $\tau$ on which 
$f$ is strictly majorized by some $\widehat{g_k}$. By
compactness of $TA$ and the fact that $(g_n) _{n=1}^\infty$ is increasing,
$f$ is majorized by some $\widehat{g_n}$ on all of $TA$.
This shows $f
\ll g$ in $\text{Cu}(A)$.

To check that $\text{LAff}_+(TA)$ is a subobject, it remains to
note that $f \ll g$ holds in $\text{LAff}_+(TA)$ if it holds in
$\text{Cu}(A)$. This follows immediately from the fact that
increasing sequential suprema in $\text{LAff}_+(TA)$ are the same
in $\text{Cu} (A)$ (see above).
\end{proof}
\end{remark}

Examples of Cuntz category semigroups
not explicitly involving a C$^*$-algebra
include $\overline{\mathbb N} := \{0, 1, 2, \ldots, \infty\}
$ and $[0,\infty]$ with the usual
order and addition taken from $\mathbb R$.
In fact, these Cuntz objects both arise from C$^*$-algebras.
Let $\text{Lsc}(X,M)$
denote the collection of lower semicontinuous functions from a
space $X$ to a Cuntz category semigroup $M$, equipped
with pointwise order and addition. Then
$\text{Lsc}(X,M)$ is a Cuntz category semigroup whenever $X$ is a
second-countable compact Hausdorff space with finite covering
dimension and $M$ is countably based (\cite[Theorem
5.17]{aps}). It was also shown in \cite[Corollary 4.22]{vilalta}
that $\text{Lsc}(X,\overline{\mathbb N})$ is a Cuntz category object 
when $X$ is a compact metric space.
It was proved in
\cite[Theorem 6.11]{ivanescu} and \cite[Theorem 10.1]{ce} 
that $\text{Cu}(C_0(X)) \cong 
\text{Lsc}(X,\overline{\mathbb N})$, via the rank map, if $X$ is $[0,1]$
or $(0,1]$.
More generally, if $X$ is a locally compact Hausdorff space
with covering dimension at most two and $\widecheck H_2(K) = 0$
(C\v ech 
cohomology with integer coefficients) for
every compact subset $K$ of $X$, then again (via the rank map)
$\text{Cu}(C_0(X)) \cong
\text{Lsc}_\sigma(X,\overline{\mathbb N})$, the ordered semigroup of
lower semicontinuous 
extended positive integer-valued functions $f$ such that $f ^{-1}
(k,\infty]$ is $\sigma$-compact for each $k \in \mathbb N$
(\cite[Theorem 1.1]{robert_cu}).
If $C_0(X)$ is separable,
these are the exact conditions needed on $X$ for
$\text{Cu}(C_0(X))$ to be isomorphic (via the rank map) to 
$\text{Lsc}(X,\overline{\mathbb N})$ (\cite[Theorem 1.3]{robert_cu}).
Note
that $X$ is a hereditarily Lindel\"of locally compact Hausdorff
space (as is the case when $C_0(X)$ is separable) if, and only if,
$\text{Lsc}_\sigma(X,\overline{\mathbb N}) =
\text{Lsc}(X,\overline{\mathbb N})$.
In general,
$\text{Lsc}_\sigma(X,\overline{\mathbb N})$ is a Cuntz category
object, and in fact (Theorem 3 of \cite{elliott_toy}) a
subobject of $\text{Cu}(C_0(X))$. 

\begin{proposition}
  Let $X$ be a locally compact Hausdorff space. Then 
  $\text{Lsc}_\sigma(X,\overline{\mathbb N})$ is a Cuntz category
  object.
  \label{lsx}
\end{proposition}

\begin{proof}


  Let us first recall a characterization of compact containment for
  $\sigma$-compact elements in the
  lattice of open subsets of $X$, ordered by
  inclusion. If there is a compact subset $K$ of $X$ sitting in
  between open sets $U$ and $V$, then $U$ is
  compactly contained in $V$.
  For the converse, suppose $U$ is compactly contained in $V$ and
  that $V$ is $\sigma$-compact.
  Since $X$ is locally compact Hausdorff, $V$ is locally compact 
  and by assumption, $V$ is $\sigma$-compact so
  there exist an increasing
  sequence of open sets $(V_k) _{k=1}^\infty$ and an increasing
  sequence of 
  compact sets $(K_k)_{k=1}^\infty$ such that $V_k \subseteq K_k
  \subseteq V_{k+1}$ for each $k$ and $\sup_k K_k = V$. Since $U$
  is compactly contained in $V$, there exists a $k$ such that $U
  \subseteq V_k \subseteq K_k \subseteq V$. 

  In the latter construction, each $V_k$ can be chosen to be
  $\sigma$-compact. Set $V_1 := \varnothing$.
  By Urysohn's lemma, there exists a continuous
  function $f_{k+1}$ from $X$ into
  $[0,1]$  which is equal to 1 on $K_k$ and is equal to 0
  on $X \backslash V_{k+1}$. The support of $f_{k+1}$,
  with which we replace $V_{k+1}$,
  is an $F_\sigma$-set sitting in between $K_k$ and
  $K_{k+1}$, and is therefore $\sigma$-compact.

  (O1) Let $(f_n) _{n=1}^\infty$ is an increasing sequence in
  $\text{Lsc}_\sigma(X,\overline{\mathbb N})$ and denote by
  $f := \sup_n f_n$. Then
  for each $k \in \mathbb N$, $\{f > k\} =
  \bigcup_{n=1}^\infty \{f_n > k\}$. (Abbreviated notation.) Since $f_n$ is lower
  semicontinuous, each $\{f_n > k\}$ is open and so $f$ is lower
  semicontinuous.
  Since each $\{f_n > k\}$ is $\sigma$-compact, a diagonalization
  argument shows $\{f > k\}$
  is $\sigma$-compact. This shows
  $\text{Lsc}_\sigma(X,\overline{\mathbb N})$ is closed under
  suprema of increasing sequences.

  (O4) is easily verified. It is also easy to see that
  if $\chi_U$ and $\chi_V$ are characteristic functions in
  $\text{Lsc}_\sigma(X,\overline{\mathbb N})$, then $\chi_U \ll
  \chi_V$ if, and only if, $U \ll V$.
  We will now show that if $f, g \in
 \text{Lsc}_\sigma(X,\overline{\mathbb N})$, then $f \ll g$ if,
 and only if, $f$ is bounded and $\{f > k\} \ll \{g > k\}$ for
 each $k$. Suppose $f$ is compactly contained in $g$. 
Each 
 $\{g>k\}$ is the supremum of a rapidly
 increasing sequence of $\sigma$-compact open sets
 $(U_{k,i})_{i=0}^\infty$ and so there exist compact sets
 $K_{k,i}$ such that $U_{k,i} \subseteq K_{k,i} \subseteq
 U_{k+1,i}$ for each $i$. Moreover, by taking finite unions
 we may choose $U_{k,i}$ and $K_{k,i}$ so that $U_{k,i} \supseteq
 U_{k+1,i}$ and $K_{k,i} \supseteq K_{k+1,i}$ whenever $k \leq
 i$.
 \[
\begin{tikzcd}
{U_{0,0}} \arrow[r, hook] & {K_{0,0}} \arrow[r, hook] & {U_{0,1}}
\arrow[r, hook]                 & {K_{0,1}} \arrow[r, hook]
& {U_{0,2}} \arrow[r, hook]                 & {K_{0,2}}
&  \\
                          &                           & {U_{1,1}}
                          \arrow[u, hook] \arrow[r, hook] &
                          {K_{1,1}} \arrow[u, hook] \arrow[r,
                          hook] & {U_{1,2}} \arrow[u, hook]
                          \arrow[r, hook] & {K_{1,2}} \arrow[u,
                          hook] & \cdots \\
                          &                           &
                          &
                          & {U_{2,2}} \arrow[u, hook] \arrow[r,
                          hook] & {K_{2,2}} \arrow[u, hook] &
\end{tikzcd}
 \]
 The supremum of the increasing sequence of functions
 $(\sum_{k=0}^i 
 \chi_{U_{k,i}}) _{i=0}^\infty$ is $g$. Since $f$ is
 compactly contained in $g$, $f$ is majorized by $\sum_{k=0}^n
 \chi_{U_{k,n}}$ for some $n$ and so $f$ is bounded. By
 construction, $\{f > k\} \subseteq U_{k,n} \subseteq K_{k,n}
 \subseteq \{g > k\}$ for $k \leq n$ (and if $k > n$, then $\{f >
 k\} = \varnothing \ll \{g > k\}$).
This shows $\{f > k\}$ is compactly
 contained in $\{g > k\}$ for each $k$. 

Conversely, suppose $f$ is
 bounded and that $\{f > k\} \ll \{g > k\}$ for each $k$. Let an
 increasing sequence $(g_i) _{i=1}^\infty $ whose supremum
 majorizes $g$ be given. 
Since $f$ is bounded, $f
= \sum_{k=0}^n \chi_{\{f > k\}}$ for some $n$. Since $\{f > k\}
\ll \{g > k\}$
there exists a
 compact set $K_k$ which sits in between $\{f > k\}$ and $\{g >
 k\}$ for each $k$. Using that $(g_i) _{i=1}^\infty$ is increasing, and
 that $g_i$ is extended
 positive integer-valued and lower semicontinuous, by compactness
 of each $K_k$, $g_i$ eventually
 majorizes $\sum_{k=0}^n \chi_{K_k}$ and therefore $g_i$
 eventually majorizes $f$. This shows $f$ is compactly contained
 in $g$.

    (O3) Suppose $f_i \ll g_i$ for $i = 1, 2$. Then $f_1 + f_2$ is
  bounded (since $f_1$ and $f_2$ are bounded).
  There is, for each $k$, a compact
  set $K_{i,k}$ sitting in between $\{f_i > k\}$ and $\{g_i > k\}$
  for $i = 1, 2$. Taking the convention that $K_{i,-1} = X$, we
  have 
  \begin{align*}
    \{f_1 + f_2 > k\}
    & = \bigcup_{i=0}^{k+1} \{f_1 > k-i\} \cap \{f_2 > i-1\} \\
    & \subseteq \bigcup_{i=0}^{k+1} \left ( K_{1,k-i} \cap K_{2,i-1}
      \right)  \\
    & \subseteq \bigcup_{i=0}^{k+1} \{g_1 > k-i\} \cap \{g_2 >
      i-1\}
      = \{g_1 + g_2 > k\}.
  \end{align*}
  This shows  that
  $\{f_1 + f_2 > k\} \ll \{g_1 + g_2 > k\}$ for each $k$, so (see
  proof of (O4) above)
  $f_1 + f_2$ is compactly contained in $g_1 + g_2$.

  (O2) Let $g \in
  \text{Lsc}_\sigma(X,\overline{\mathbb N})$ be given. Each $\{g >
  k\}$ is the supremum of a rapidly increasing sequence of open
  sets $(U_{k,i})_{i=0}^\infty$ and so $g$ is the supremum of
  the increasing sequence 
  $(\sum_{k=0}^n \chi_{U_{k,n}})_{n=0}^\infty$. By (O3), this
  sequence is rapidly increasing.
\end{proof}

  \begin{theorem}[{\cite[Theorem 1.1]{robert_cu}}]
    Let $X$ be a locally compact Hausdorff space of covering
    dimension at most two with $\check
    H_2(K) = 0$ 
    for
    every compact subset $K$ of $X$.
    Then $\text{Cu}(C_0(X)) \cong
    \text{Lsc}_\sigma(X,\overline{\mathbb N})$, the isomorphism
    consisting of the rank map.
  \label{robert_cu}
\end{theorem}

In order to lift Cuntz category morphisms to $*$-homomorphisms, we
will need the following classification theorem of
\cite{ce}. Further developments in this direction
include \cite[Theorem 2]{robert_santiago}, \cite[Theorem 1]{ces},
and \cite[Theorem 1.0.1]{robert}. 

\begin{theorem}[{\cite[Theorem 4.1 and remark page 29]{ce}}]
  \label{classify_homs}
    Let $B$ be a stable rank one C$^*$-algebra and let $s_B$ be a
    strictly positive element of $B$. If $\sigma :
    \text{Cu}(C_0(0,1]) \to \text{Cu}(B)$ is a Cuntz category
    morphism taking $[t]$ into $[s_B]$, or into $g \leq [s_B]$,
    then there is a
    $*$-homomorphism $\pi : C_0(0,1] \to B$, which is unique up to
    approximate unitary equivalence, such that $\text{Cu}(\pi) =
    \sigma$. 
  \end{theorem}

\subsection{Embedding theorems} 

\begin{theorem}[{\cite[Theorem 2.8]{kp}}]
  \label{o2-embedding}
  If $A$ is a unital, separable, and exact C$^*$-algebra,
  then there exists a unital embedding
  of $A$ into $\mathcal O_2$.
\end{theorem}

\begin{theorem}
  \label{AF-embedding}
  If $A$ is a unital, finite, $\mathcal Z$-stable, simple,
  separable amenable C$^*$-algebra satisfying the UCT,
  then there is a
  unital embedding of $A$ into a unital, simple, separable AF
  algebra $D$ giving rise to an isomorphism of tracial simplices.
\end{theorem}

\begin{proof} Let $\mathcal Q$ denote the universal UHF algebra.
  Since there is an embedding of $A$ into $A
  \otimes \mathcal Q$ which induces an isomorphism of tracial
  simplices,
  we may suppose, without loss of generality,
  that $A$ is $\mathcal Q$-stable. (This will be used in applying
  \cite{gln} below.)
  
  Let $\rho : K_0(A) \to \text{Aff}(T(A))$
  denote the map associated to the
  canonical pairing of $K_0(A)$ with $T(A)$.
  Choose a countable dense subgroup $G$ of $\text{Aff}(T(A))$
  containing the image of $\rho$.
  By \cite{ando} and \cite{lindenstrauss}, $\text{Aff}(T(A))$ with
  the strict pontwise order has the Riesz interpolation property,
  and therefore $G$ does also.
  By
  \cite{ehs}, there exists a unital, simple, separable
  AF algebra $D$ with $K_0(D) = G$. Since $K_0(D)$ is dense in
  $\text{Aff}(T(A))$, the map $T(A) \to S(K_0(D))$ is an
  isomorphism. Since $D$ is an AF algebra, $T(D)  = S(K_0(D))$, and so
  we have an isomorphism $\Phi : T(D) \to T(A)$. We now have
  compatible maps $\rho: K_0(A) \to K_0(D)$ and $\Phi : T(D) \to
  T(A)$. Furthermore, $A$ and we may suppose also $D$ are
  $\mathcal Q$-stable.
  By \cite[Theorem 1.1]{egln}, \cite[Theorem A]{tww}, and
  \cite[Theorem A]{cetww}, $A \otimes \mathcal Q$
  has generalized tracial rank at most one -- and therefore, $A$
  does too as it is $\mathcal Q$-stable.
  The only additional constituent of the invariant
  (\cite[Definition 2.4]{gln}) used in the homomorphism theorem
  \cite[Corollary 21.11]{gln} is a $K_1$-map which we can take to
  be zero. The corresponding algebra map is the desired embedding.
  \end{proof}

\subsection{Uniqueness theorems} \label{uniqueness}

A couple of the key technical ingredients used
in the proof of Theorem
\ref{thm} are the following
uniqueness theorems for maps from unital, separable amenable
C$^*$-algebras into unital, simple, separable, $\mathcal
Z$-stable C$^*$-algebras. These results were developed in
\cite{bbstww} with certain restrictions on the trace
space. These tracial assumptions were removed in
\cite{cetww}. For the statement involving a Kirchberg algebra
(i.e.~ a unital, purely
  infinite, simple, separable amenable C$^*$-algebra) as
the codomain, we specialize to the case that the maps are
injective $*$-homomorphisms.
Such maps $\varphi$
induce injective c.p.c.~ order zero maps $(\varphi-t)_+$ for each
$t \in [0,1)$ and so
\cite[Corollary 9.11]{bbstww} is applicable to them. (To see
this, it suffices to show that $(h_\varphi-t)_+$ is a non-zero
scalar multiple of the projection $\pi(1_A)$ where notation is as in Theorem
\ref{structure} (note that, here, $\pi_\varphi = \varphi$). Since
$h_\varphi$ is a projection (namely, $\pi_\varphi(1_A)$),
$(h_\varphi -t)_+  = (1-t) h_\varphi$.)

\begin{theorem}[{\cite[Corollary 9.11]{bbstww}}]
  Let $A$ be a unital, separable, amenable C$^*$-algebra and let
  $B$ be a Kirchberg algebra. Let $\phi_1, \phi_2 : A \to
  B_\omega$ be a pair of injective $*$-homomorphisms.   Then,
  there exist contractions $v_i, w_i \in
  B_\omega$, $i = 1, 2$, such that
  \begin{align*}
    \phi_1(a) & = w_1 \phi_2(a) w_1^* + w_2 \phi_2(a) w_2^*, \\
    \phi_2(a) & = v_1 \phi_1(a) v_1^* + v_2 \phi_1(a) v_2^* 
  \end{align*}
  for all $a \in A$, 
  and
  \[
    v_1^*v_1 + v_2^*v_2 = w_1^*w_1 + w_2^*w_2 =
    1_{B_\omega}
  \]
  with $w_i^*w_i \in \phi_2(A)'$ and $v_i^*v_i
  \in \phi_1(A)'$, $i=1,2$.
  \label{kirchberg_colored}
\end{theorem}

\begin{definition} \label{coloured_equivalent}
  Let $A$ and $B$ be unital C$^*$-algebras.
  Any two order zero maps $\phi_1, \phi_2 : A \to B_\omega$
  satisfying the conclusion of Theorem
\ref{kirchberg_colored} will be said to be {\it coloured
  equivalent} (\cite[Section 6]{bbstww}).
\end{definition}

For the statement involving a finite algebra as the codomain, only
one map is assumed to be a homomorphism.

\begin{theorem}[{\cite[Theorem 6.6]{bbstww}}, \cite{cetww}]
  Let $A$ be a unital, separable 
  amenable C$^*$-algebra
  and let $B$ be a unital,
  finite, $\mathcal Z$-stable, simple,
  separable amenable C$^*$-algebra.
  Let $\phi_1 : A \to B_\omega$ be a
  totally full $*$-homomorphism (i.e., $\phi_1(a)$ generates
  $B_\omega$ as a closed two-sided ideal for every non-zero $a \in A$) and $\phi_2 : A \to B_\omega$ be
  a c.p.c.~order zero map with
  \[\tau \circ \phi_1 = \tau \circ
    \phi_2^m
  \]
  for all $\tau \in T(B_\omega)$ and all $m \in
  \mathbb N$. Then, $\phi_1$ and $\phi_2$ are coloured equivalent.
   \label{coloredclassification}
\end{theorem}

\begin{proof} By
  \cite[Theorem I]{cetww}, $B$ has
  complemented partitions of unity.
  Now the proof is exactly as in \cite[Theorem 6.6]{bbstww} using
  \cite[Lemma 4.8]{cetww} in place of \cite[Theorem
  5.5]{bbstww}. The same lemma also permits relaxing the tracial
  extreme boundary hypothesis in \cite[Theorem 6.2]{bbstww}.
\end{proof}

In practice, the $*$-homomorphism typically used in Theorem 
\ref{coloredclassification} is essentially the canonical
embedding $\iota_A: A \to A_\omega$ (e.g.~ \cite[Corollary
6.5]{bbstww} and \cite[Theorem
7.5]{bbstww}), which is totally full if $A$
is simple. To see this, let a non-zero element $a \in A$ be given.
If $A$ is unital, $\iota_A(a)$ can be cut down to any
coordinate with a projection. In the non-unital case, the
coordinate projection can be replaced with an approximate unit in a
coordinate. Since $A$ is simple, a copy of $A$ is generated in
each coordinate of $A_\omega$. This shows $\iota_A$ is totally
full. 

\section{Functionals preserving compact containment}
\label{functionals}

The main result of this section characterizes when a functional on
the Cuntz semigroup $\text{Cu}(A)$ 
arising from a faithful densely defined lower semicontinuous
traces on a commutative
C$^*$-algebra, $A$, 
with the property that the rank map gives rise to an isomorphism
$\text{Cu}(A)
\cong \text{Lsc}_\sigma(\widehat A,\overline{\mathbb N})$,
preserves compact
containment. Let us first establish some necessary conditions for a
functional to preserve compact containment in a
more general context.
Suppose $\lambda$ is a functional (i.e.~a generalized Cuntz
category morphism into $[0,\infty]$) on a positive Cuntz category
object $M$. 
In order for $\lambda$ to preserve compact
containment, 
$\lambda$ must be finite on any
element that is compactly contained in some element of
$M$. Moreover, $\lambda$ must be finite on any element that is
majorized by a finite sum of such elements.
It is also necessary that 
$\lambda$ vanish on compact elements.
Now recall, for instance from \cite[Proposition 4.2]{ers},
that a lower semicontinuous trace $\tau$ on a C$^*$-algebra
induces a functional $d_\tau$ on $\text{Cu}(A)$
(the trace of the range projection of a positive element).
(In later
sections, we will not make a distinction between $\tau$ and
$d_\tau$.)
The finiteness condition requires $d_\tau$ to be finite-valued on
the Pedersen ideal of $A \otimes \mathcal K$. 
If $\tau$ is faithful,
the latter condition requires $\text{Cu}(A)$ to have no non-zero compact
elements. If $A$ is stably finite, then, by
\cite[Corollary 3.6]{bc}, this is equivalent to $A$ being stably
projectionless. If $A$ is, moreover, commutative, this is
equivalent to the spectrum of $A$ not containing any non-empty
compact open subsets. To see this, suppose $K$ is a non-empty
compact open subset of the spectrum of $A$. Then $[\chi_K \otimes p]$ is a
non-zero compact element of $\text{Cu}(A)$ for any non-zero
projection $p \in \mathcal K$. Conversely, suppose $A
\otimes \mathcal K$ contains a non-zero projection $p$.
Then the map
$\eta : \widehat A \ni x \mapsto \|p(x)\| \in \mathbb C$
is continuous.
Since $p$
vanishes at infinity and is a non-zero projection, the same is
true of $\eta$. This shows $A$ contains a non-zero projection.

\begin{theorem}
  Suppose that $\mu$ is a faithful densely defined lower
  semicontinuous trace on 
  a commutative C$^*$-algebra, $C_0(X)$, with the property that
  $\text{Cu}(C_0(X)) \cong \text{Lsc}_\sigma(X,\overline{\mathbb
    N})$ (the isomorphism consisting of the rank map) (for
  example, $X$ is as in Theorem \ref{robert_cu}). 
  Then the functional $d_\mu: \text{Cu}(C_0(X)) \to [0,\infty]$
  preserves 
  compact containment if, and only if, $X$
  does not contain any
  non-empty compact open sets.
  \label{faithful_commutative}
\end{theorem}

\begin{proof}
 Denote by $\iota$ the natural set-theoretical mapping of
 $\text{Lsc}_\sigma(X,\overline{\mathbb N})$ into
 $\text{Cu}(C_0(X))$:
 \[
   f = \sum_{k=0}^\infty \chi_{\{f > k\}} \mapsto \sum_{k=0}^\infty
   [f_k \otimes p_k],
 \]
 where $f_k$ is a positive function in $C_0(X)$ with support
 equal to $\{f > k\}$ and $(p_k) _{k=0}^\infty$ is a family of
 mutually orthogonal rank-one projections in $\mathcal K$.
 By hypothesis, $\text{rank} : \text{Cu}(C_0(X)) \to
 \text{Lsc}_\sigma(X,\overline{\mathbb N})$ is an
 isomorphism. It follows, since as easily checked,
 $\text{rank} \circ \iota =
 \text{id}_{\text{Lsc}_\sigma(X,\overline{\mathbb N})}$, 
 that $\iota$ is an isomorphism.
 Denote by $\mu :
 \text{Lsc}_\sigma(X,\overline{\mathbb N}) \to [0,\infty]$ the
 functional $d_\mu \circ \iota$
 and denote also by $\mu$ the
 locally finite (as the Pedersen ideal of $C_0(X)$ is $C_c(X)$
 and by \cite[Theorem 1.3]{pedersen})
 extended Radon measure
 induced by the densely defined lower semicontinuous
 trace $\mu$.
 Since $d_\mu = \mu \circ \operatorname{rank}$,
 to show that $d_\mu$ preserves
 compact containment, it suffices to show that $\mu$ preserves
 compact containment. Suppose $f$ and $g$ are elements of
 $\text{Lsc}_\sigma(X,\overline{\mathbb N})$ and that $f \ll
 g$. If $\mu(f) = 0$, then $\mu(f) \ll \mu(g)$ because $0$ is compactly
 contained in every element of $[0,\infty]$. Let us now consider the
 case where $\mu(f) > 0$.
 Since $f$ is compactly contained in $g$,
 $f = \sum_{k=0}^n
 \chi_{\{f > k\}}$ for some $n$
 and
 there exists a compact set $K \subseteq X$ such that (with
 abbreviated notation)
 $\{f > 0\} \subseteq K \subseteq \{g > 0\}$, and so $\mu(f) \leq
 (n+1)\mu(K)< \infty$.
 Now suppose $\mu(f) = \mu(g)$, i.e.,
  $\sum_{k=0}^\infty \mu( \{f > k\}) =
 \sum_{k=0}^\infty \mu( \{g > k\})$. Since $f \leq g$ and since
 $\mu(f)$ is finite, we
 have $\mu(\{f > k\}) = \mu(\{g >k\})$ for each $k$. 
 If $K = \{g >
 0\}$, then since $X$ contains no nonempty compact open sets, $\{g
 > 0\}$ is empty. This implies $g$ is zero and hence $f$ is zero,
 which is a
 contradiction.
 In other words, $K$ is properly contained in
 $\{g > 0\}$, and so there is a
nonempty open subset of $\{g > 0\}$ which is disjoint
 from $K$ and therefore disjoint from $\{f > 0\}$ as well.
 Since $\mu$ is faithful,
 $\mu(\{f > 0\}) < \mu(\{g > 0\})$ and so (since $f \leq g$ and
 $\mu(f)$ is finite) $\mu(f) < \mu(g)$.
\end{proof}

\begin{remark} \label{aps_proof}
  A different proof of Theorem
  \ref{faithful_commutative},
  in the case of the half-open interval (the case pertinent to this paper),
  using the characterization of compact containment given in 
  \cite{aps} is possible.

  \begin{proof} 
    Suppose $f, g \in
  \text{Lsc}((0,1],\overline{\mathbb N})$ are such that $f \ll
  g$.
  Let $\mu$ be as in the proof of Theorem
  \ref{faithful_commutative} and 
  suppose $f$ is non-zero. Extend both $f$ and $g$ to be zero at
  zero and call these extensions $\widetilde f$ and $\widetilde
  g$. It is easily checked that these extensions are lower
  semicontinuous and that $\widetilde f$ remains compactly
  contained in $\widetilde g$ in
  $\text{Lsc}([0,1],\overline{\mathbb N})$ (one could also use
  \cite[Theorem 6.11]{ivanescu}, \cite[Theorem 10.1]{ce},
  and \cite[Theorem 5]{crs}). 
  An application of \cite[Proposition 5.5]{aps} at
  the point zero shows that $\widetilde f$ is zero on a
  neighbourhood of zero because $\widetilde g(0) = 0$.
  This shows $f$ is compactly supported, and since $f$ is bounded,
  $\mu(f)$ is finite.
  Since $f$ is
  non-zero and extended positive integer-valued, its extension
  $\widetilde f$ is necessarily discontinuous. Its first point of
  discontinuity, $x_0$, 
  must occur within $(0,1)$. By
  \cite[Proposition 5.5]{aps} applied to the point $x_0$, there
  exists a neighbourhood $U$ of $x_0$ which is contained in
  $(0,1)$ and a non-zero constant $c
  \in 
  \overline{\mathbb N}$ such that $\widetilde f \leq c \ll
  \widetilde g$ on $U$. Since $\widetilde f$ is zero to the
  left of $x_0$, and since $\widetilde f \leq \widetilde g$, the
  integral of $f$ on $U$ is strictly less than the
  integral of $c$ on $U$, and hence also of $g$. This implies
  $\mu(f) < \mu(g)$.
\end{proof}
\end{remark} 

\begin{remark} \label{faithful}
  The characterization of compact containment in \cite{cei} can be
  used to show that if $\tau$ is a faithful densely defined lower
  semicontinuous trace on a 
  C$^*$-algebra $A$, then $d_\tau$
  preserves compact containment if, and only if, $\text{Cu}(A)$
  does not contain non-zero compact elements.
  Since the only case of interest in this paper is covered by
  Theorem \ref{faithful_commutative}, we will not prove this.
\end{remark}

\begin{remark}
  A simple reduction to the faithful case shows
  that if $\tau$ is a (not necessarily faithful) densely defined lower semicontinuous trace
  on $A$ and $\text{Cu}(A/N_\tau)$ does not contain non-zero
  compact elements, where $N_\tau$ denotes the kernel of $\tau$,
  then $d_\tau$ preserves compact containment.
  
  \begin{proof} 
  Denote by $\pi_{N_\tau} : A \to A/N_\tau$ the canonical quotient
  map and let $\eta : (A/N_\tau)_+ \to [0,\infty]$ be the faithful
  densely defined lower semicontinuous 
  trace induced by $\tau$. Suppose $\text{Cu}(A/N_\tau)$ contains
  no non-zero compact elements. Then, by Remark \ref{faithful}, the functional $d_\eta : \text{Cu}(A/N_\tau)
  \to [0,\infty]$ preserves compact containment and
  by
  \cite[Theorem 2]{cei},
  $\pi_{N_\tau}$ induces a Cuntz category morphism
  $\text{Cu}(\pi_{N_{\tau}}) : \text{Cu}(A) \to
  \text{Cu}(A/N_\tau)$. Since 
    the functional $d_\tau$ factors
  through $\text{Cu}(A/N_\tau)$, via the commutative diagram
  \[
    \begin{tikzcd}
      \text{Cu}(A) \arrow[rr, "\textstyle d_\tau"] \arrow[rd, "\textstyle \text{Cu}(\pi_{N_\tau})"'] &                                                      & {[0,\infty]} \\
      & \text{Cu}(A/N_\tau) \arrow[ru, "\textstyle d_\eta"'] &             
    \end{tikzcd},
  \]
it preserves compact
  containment. 
\end{proof}

\end{remark}

\begin{remark}
  It is possible for the functional $d_\tau$ induced
  by a lower semicontinuous trace $\tau$ on a C$^*$-algebra $A$
  to fail to preserve
  compact containment
  even if
  $d_\tau$ is assumed to vanish on every compact element of
  $\text{Cu}(A)$.
  For example, let $\mu$ be the trace
  induced by the measure
  on $(0,1]$ which is the Lebesgue
  measure on $(1/2,1]$ and zero on $(0,1/2]$.
  The condition that $d_\mu$ vanishes on every compact element of
  $\text{Cu}(C_0(0,1])$ is automatic since $C_0(0,1]$ is stably
  projectionless.  

To show that
  $d_\mu$ fails to preserve compact containment, it suffices, by
  the proof of Theorem \ref{faithful_commutative}, to show that $\mu$ fails to
  preserve 
  compact containment.
  By a characterization of compact containment in the
  proof of
  Proposition \ref{lsx},
  $\chi_{(1/2,1]}$ is compactly contained in $\chi_{(0,1]}$, but
  $\mu((1/2,1]) = \mu((0,1])$.
\end{remark}

\section{Coloured Classification of C$^*$-algebras} \label{main_section}

\begin{theorem} \label{prop} 
  Any two unital C$^*$-algebras which are coloured
  isomorphic have isomorphic tracial simplices. 
  More precisely, the order
  zero maps implementing a coloured isomorphism
  induce mutually inverse isomorphisms of tracial
  simplices (as described in Section \ref{traces}).
  \label{preserve_traces}
\end{theorem}

\begin{proof}
  Suppose $A$ and $B$ are unital C$^*$-algebras which are coloured
  isomorphic with notation as in \ref{colored}. Since $\varphi : A_\omega \to
  B_\omega$ and $\psi : B_\omega \to A_\omega$ are c.p.c.~ order
  zero maps, they induce mappings of bounded traces $\varphi^* :
  \mathbb R_+ T(B_\omega) \to \mathbb R_+(A_\omega)$ and $\psi^* :
  \mathbb R_+ T(A_\omega) \to \mathbb R_+(B_\omega)$
  (\cite[Corollary 3.4]{wz_cpc}). In fact, $\varphi^*$ and
  $\psi^*$ are
  mutually inverse affine isomorphisms of 
  the cones of bounded traces. To see this, 
  let $\tau \in T(A_\omega)$ and $a \in A_\omega$ be
  given. Then
  \begin{align*}
    \tau(a)
    & \stackrel{(\ref{coloredid})}{=}
       \tau \left ( \sum_{i=1}^m u_i \psi \varphi(a) u_i^* \right)
      = \sum_{i=1}^m \tau(u_i \psi \varphi (a) u_i^*) \\
    & = \sum_{i=1}^m \tau(\psi \varphi (a) u_i^* u_i)
      \stackrel{(\ref{coloredpartition})}{=} 
      \tau(\psi \varphi (a)) = \varphi^* \psi^*(\tau)(a).
  \end{align*}
  A symmetric calculation shows $\psi^* \varphi^* (\tau) =\tau$
  for each $\tau \in T(B_\omega)$.

  Let us now show that the cones of bounded traces on $A$ and $B$
  are affinely isomorphic.
  Denote by $c_B : \mathbb R_+T(B) \to \mathbb
  R_+T(B_\omega)$ the embedding of bounded traces on $B$
  as bounded limit traces on $B_\omega$ 
  and by $\iota_A : A \to A_\omega$ the canonical 
  embedding of $A$ into $A_\omega$. 
  Denote by $\iota_A^*$ the dual map induced on bounded traces and define $\Phi :
  \mathbb R_+T(B) \to \mathbb R_+T(A)$ to be the composition
  \[ 
    \begin{tikzcd}
\mathbb R_+T(B) \arrow[r, "\textstyle c_B"] & \mathbb R_+T(B_\omega)
\arrow[r, "\textstyle \varphi^*"] & \mathbb R_+T(A_\omega) \arrow[r,
"\textstyle \iota_A^*"] & \mathbb R_+T(A),
\end{tikzcd}
\]
and define $c_A : \mathbb R_+T(A) \to \mathbb R_+T(A_\omega)$, $\iota_B: B \to
B_\omega$, and $\Psi : \mathbb R_+T(A) \to \mathbb R_+T(B)$ similarly.
Then $\Phi$ and $\Psi$ are weak* continuous and affine, and
for $\tau \in
\mathbb R_+T(B)$ and $b \in B$,
\begin{align*}
  \Psi \Phi (\tau)(b)
  & = \iota_B^* \psi^* c_A \iota_A^* \varphi^* c_B (\tau) \\
  & = \iota_B^* \psi^* \varphi^* c_B(\tau) = \iota_B^* c_B(\tau) = \tau.
\end{align*} 
The second equality follows from the assumption that 
$\varphi^*$ takes constant limit traces to constant limit traces
so that
$c_A \iota_A^* \varphi^* c_B (\tau) = \varphi^*
c_B(\tau)$. The third equality follows from the previous
paragraph.
A symmetric calculation shows that $\Phi \Psi(\tau) = \tau$ for
$\tau \in \mathbb R_+T(A)$.
This shows $\Phi$ and $\Psi$ are
mutually inverse isomorphisms of the topological convex cones
$\mathbb R_+ T(A)$ and $\mathbb R_+ T(B)$.

To see that $\Phi$ and $\Psi$ preserve the tracial simplices, i.e., are
isometries,
note that they are contractions since they are the compositions of
contractions. By the preceding paragraph, $\Phi \Psi =
\text{id}_{T(A)}$ and so we have $1 = \|\Phi
\Psi\| \leq \|\Phi\| \|\Psi\| \leq 1$ which implies that the norms
of $\Phi$ and $\Psi$ are both one so these maps are isometries
and therefore constitute an
isomorphism of the tracial simplices. \end{proof}

\begin{remark} Coloured isomorphism
  (as defined in Section \ref{colored}), extended in a natural way
  to the non-unital case,
  also preserves, up to
  isometric isomorphism, 
  the topological convex cone of lower semicontinuous traces of
  \cite{ers}.
  To see this, we will need that c.p.c.~ order zero
  maps induce mappings of lower semicontinuous traces.
  This follows from \cite[Corollary
  3.4]{wz_cpc}. The result now follows by
  extending the maps considered in the proof of 
  Theorem \ref{preserve_traces} to the cone of lower
  semicontinuous traces rather than just the cone of bounded
  traces.

  From this one can deduce that coloured isomorphism, in the
  present sense, preserves
  ideal lattices, up to isomorphism
  (cf.~ \cite[Theorem 5.4.9]{castillejos}). 
  Let notation be as in the proof of Theorem
  \ref{preserve_traces} with $\Psi$ assumed to be an isomorphism
  of cones of lower semicontinuous traces.
  Recall, from \cite{ers}, that there
  is an order-reversing bijection $\alpha_A$
  between closed two-sided ideals
  $I$ of $A$ and lower semicontinuous traces on $A$ taking
  only the values 0 and $\infty$ given by
  \[
    I \mapsto \tau_I(x) :=     \begin{cases}
      0 & \text{if}   \hspace{5 pt}    x \in I^+\\
      \infty & \text{if}    \hspace{5 pt}   x \notin I^+
    \end{cases}.
  \]
  Such traces are characterized by the property $\tau + \tau =
  \tau$, and so $\Psi$ restricts to a mapping between these
  traces. Therefore, the
  composition $\alpha_B ^{-1} \circ \Psi \circ \alpha_A$ gives an
  order-preserving bijection of ideal lattices. If the prime
  ideals of $A$ and $B$ are primitive, as is the case with separable
  or postliminary C$^*$-algebras (\cite[Theorem A.50]{rw}, 
  \cite[Theorem 4.4.5]{dixmier}, and \cite[Theorem 4.3.5]{dixmier}),
  then the primitive ideal spaces of $A$ and $B$ are
  homeomorphic. By Gelfand duality,
  it follows that coloured isomorphism (in the present sense)
  is a rigid notion
  for commutative C$^*$-algebras (i.e.~ if $A$ is a commutative
  C$^*$-algebra which is coloured isomorphic, in the present
  sense, to a C$^*$-algebra
  $B$, then $A$ is isomorphic to $B$)
  (cf. \cite[Proposition
  5.4.13]{castillejos}).
  It is easily seen that 
  coloured isomorphism, in the present sense,
  is also a rigid notion for
  finite-dimensional C$^*$-algebras 
  (cf.~ \cite[Proposition 5.4.12]{castillejos}).
\end{remark}

\begin{theorem}   \label{induce_tensor_trace}
  Let $B$ be a unital, simple, separable, exact,
  $\mathcal Z$-stable C$^*$-algebra
  with stable rank one and let $D$ be a simple, unital AF algebra.
  Let $\mu$ be a faithful trace on $C_0(0,1]$ with norm at most
  one and 
  let 
  $\Phi : T(B) \to T(D)$ be a  continuous
  affine map of simplices. There exists a c.p.c.~ order zero map
  $\varphi : D \to B$
  which satisfies the identity
  \begin{align}
    \tau \varphi^n = \mu(t^n) \Phi(\tau)     \label{tensor_trace}
  \end{align}
  for all $n \in \mathbb N$ and for all $\tau \in T(B)$.

\end{theorem}

\begin{proof}
    
  Let $(D_i)_{i=1}^\infty$ be an increasing sequence of
  finite-dimensional C$^*$-subalgebras of $D$ such that
  $\bigcup_{i=1}^\infty D_i$ is dense in $D$. Denote by $\iota_i$
  the embedding of $D_i$ into $D$.
  We claim that the map
  \[
    \sigma_i :\text{Cu}(C_0(0,1] \otimes D_i) \to
    \text{LAff}_+( TB)
  \]
  determined by the rule
  \[
    [d] \mapsto (\mu \otimes (\iota_i^* \circ \Phi(\cdot)))[d]
  \]
  is a Cuntz category morphism.

  We first show that 
  $\sigma_i$ in fact maps into $\text{LAff}_+(TB)$. 
  The image of $[d] \in \text{Cu} (C_0(0,1] \otimes D_i)$
  is clearly a positive real-valued affine
  function on $TB$.
  By Section 5 of \cite{ers}, $\sigma_i([d])$ is lower
  semicontinuous.
  By assumption, $\mu$ is faithful and $\iota_i^* \circ
  \Phi(\tau)$, being the restriction of a faithful trace on $D$
  (as $D$ is simple), is also faithful for each $\tau \in
  TD_i$. Therefore, $(\mu \otimes (\iota_i^* \circ
  \Phi(\cdot)))[d]$ is (pointwise)
  strictly positive whenever $[d]$ is
  non-zero.  
  Since $B$ is unital and
  separable, $TB$ is metrizable. Therefore, by \cite[Corollary
  I.1.4]{alfsen} and \cite[Lemma 3.6]{thiel}, $\sigma_i([d])$ is the
  pointwise supremum of an increasing sequence of continuous
  finite-valued functions in $\text{LAff}_+(TB)$. This shows that
  $\sigma_i([d])$ is an element of $\text{LAff}_+(TB)$.

  That $\sigma_i$ is a generalized Cuntz category
  morphism
  follows from Section 4 of \cite{ers}, so all that remains is
  to show $\sigma_i$ preserves compact containment. 
  Since $D_i$ is
  a finite-dimensional C$^*$-algebra, we may identify it with a
  finite direct sum of matrix  algebras,
  $\bigoplus_{j=1}^k M_{n_j}$. 
Using that $C_0(0,1]
  \otimes (\bigoplus_{j=1}^k M_{n_j})$ is isomorphic to
  $\bigoplus_{j=1}^k 
  (C_0(0,1] \otimes M_{n_j})$, we make the
  identification $\text{Cu}(C_0(0,1] \otimes D_i) =
  \bigoplus_{j=1}^k \text{Cu}(C_0(0,1] \otimes M_{n_j})$. 
  By \cite[Appendix 6]{cei}, the embedding of
  $C_0(0,1]$ into the upper-left corner of $C_0(0,1] \otimes
  M_{n_j}$ induces an isomorphism at the level of the Cuntz
  category.   For each $j$,
  we denote by $e_j$ the non-zero minimal projection in
  the upper-left corner of $M_{n_j}$, and we denote by $\rho_j$
  the pure state on $M_{n_j}$. 
  These maps induce an isomorphism
  $\bigoplus_{j=1}^k \text{Cu}(C_0(0,1] \otimes M_{n_j})
  \cong
  \bigoplus_{j=1}^k \text{Cu} (C_0(0,1])$. By Theorem
  \ref{robert_cu}, $\text{Cu}(C_0(0,1]) \cong
  \text{Lsc}((0,1],\overline{\mathbb N})$ and so under these
  isomorphisms, an arbitrary element $[d] \in
  \text{Cu}(C_0(0,1] \otimes D_i)$
  is a direct sum of elements $d_1,
  \ldots, d_k$ of $\text{Lsc}((0,1],\overline{\mathbb N})$.
  Pick $d_j^i \in C_0(0,1]$ to have support exactly $\{d_j > i\}$
  ($j = 1, \ldots, k$ and $i \in \mathbb N$).
  Let $(p_i)_{i=0}^\infty$ be a sequence of mutually orthogonal
  rank-one projections and let $\iota:
  \text{Lsc}((0,1],\overline{\mathbb N}) \to \text{Cu} (C_0(0,1])$
  and $\mu : \text{Lsc}((0,1],\overline{\mathbb N}) \to
  [0,\infty]$ be as in the proof of Theorem
  \ref{faithful_commutative}. The direct sum of copies of
  $\iota$ gives an isomorphism 
  $\bigoplus_{j=1}^k \text{Lsc}((0,1],\overline{\mathbb N}) \to
  \bigoplus_{j=1}^k\text{Cu}(C_0(0,1])$. Under the aforementioned
  isomorphisms, we have  
  \[
    [d] = \sum_{j=1}^k \sum_{i=0}^\infty [d_j^i \otimes e_j
    \otimes p_i].
  \]
  Therefore,
  \begin{align*} 
    \sigma_i([d])
    & = \sum_{j=1}^k \sum_{i=0}^\infty \mu(d_j^i) (\iota_i^* \circ
      \Phi(\cdot))(e_j)
      \text{Tr} (p_i) \\
    & = \sum_{j=1}^k \sum_{i=0}^\infty \mu(d_j^i) (\iota_i^* \circ
      \Phi(\cdot))(e_j)       
    \\
    & = \sum_{j=1}^k \mu(d_j) (\iota_i^* \circ \Phi(\cdot))(e_j).
  \end{align*}

  It is readily seen from this calculation that $\sigma_i([d])$ is
  continuous. If $\sigma_i([d])$ is bounded,
  then $\sigma_i([d])$ is an
  affine extension of the continuous function taking on the real
  values $\mu(d_j)$ at the extreme points $\rho_j$ of
  $TD_i$. Since $D_i$ is finite-dimensional, $TD_i$ is a Bauer
  simplex and so $\sigma_i([d])$ is continuous. If $\sigma_i([d])$
  is not bounded, then $\mu(d_j)$ must be equal to $\infty$ for
  some $j$. Since $\iota_i^* \circ \Phi(\cdot)$ maps into the
  faithful part of $TD_i$, $\sigma_i([d])$ is the continuous
  function which is constant and equal to $\infty$. 
  
  Suppose $[f] \ll [g]$ in $\text{Cu}(C_0(0,1] \otimes D_i)$. Let
  notation be as in the preceding two paragraphs for both $[f]$ and
  $[g]$. If $[f]$ is equal to zero, then $\sigma_i([f])$ is also
  equal to zero and it is compactly contained in every element of
  $\text{LAff}_+( TB)$ (and in particular, $\sigma_i([g])$).
  If $[f]$ is non-zero, then $f_j \ll g_j$ for each $j$ and
  $f_j$ must be nonzero for some $j$. For this particular $j$,
  $\mu(f_j) < \mu(g_j)$, by faithfulness of $\mu$, and for each
  $j$, $\mu(f_j)$ is 
  finite and $\mu(f_j) \leq \mu(g_j)$ (these facts are
  established in
  the proof of Theorem \ref{faithful_commutative}).
  It follows from the above computation and the fact that
  $\iota_i^* \circ \Phi(\tau) > 0$ ($\tau \in TB$) that 
  $\sigma_i([f]) < \sigma_i([g])$. By the previous
  paragraph, this also shows that $\sigma_i([f])$ is continuous and
  finite-valued. Since lower semicontinuous functions attain their
  infima on compact sets, it follows from the characterization of
  compact containment in Section 2 of \cite{ers} (or \cite[Lemma
  3.6]{thiel}) that $\sigma_i([f])$ is compactly contained in
  $\sigma_i([g])$. This shows $\sigma_i$ is a Cuntz category
  morphism.

  Since $\text{LAff}_+(TB)$ is (in a natural way)
  a subobject of $\text{Cu}(B)$
  (Remark \ref{subobject}), $\sigma_i$ extends to a Cuntz category
  morphism into $\text{Cu}(B)$ which we denote again by
  $\sigma_i$. In order to lift $\sigma_i$ to a $*$-homomorphism, 
  let us check
  that $\sigma_i$ takes a strictly positive element into an
  element majorized by a strictly positive element:
  \begin{align*}
    \sigma_i([t \otimes 1_{D_i}])
    & = \sum_{j=1}^k \mu(\chi_{(0,1]})(\iota_i^* \circ
      \Phi(\cdot))[1_{M_{n_j}}]
      \leq \sum_{j=1}^k (\iota_i^* \circ \Phi(\cdot))[1_{M_{n_j}}] \\
    & = \sum_{j=1}^k \Phi(\cdot)[\iota_i(1_{M_{n_j}})]
      = \Phi(\cdot) [\iota_i(1_{D_i})]
      \leq \widehat{[1_B]}.
  \end{align*}
  In the first inequality, we have used that $\mu$ is of norm at
  most one.
  By Theorem \ref{classify_homs}, i.e.,
  by \cite[Theorem 4.1]{ce}, as modified in the remark on
  page 29 of \cite{ce} (to replace $[0,1]$ with $(0,1]$), for each
  $i$ there exists a
  $*$-homomorphism $\pi_i : C_0(0,1] \otimes D_i \to B$, which is
  unique up to approximate unitary equivalence, such that
  $\text{Cu} (\pi_i) = \sigma_i$.
  Let us now denote by $\iota_i$ the embedding of
  $C_0(0,1] \otimes D_i$ into $C_0(0,1] \otimes D_{i+1}$.
  Then with the embeddings $\text{Cu}(\iota_i)$ on the top row and
  the identity map $\text{Cu}(\text{id}_B)$ on the bottom row,
  we have the following one-sided intertwining at the level
  of the Cuntz category:
  \[
    \begin{tikzcd}
      {\text{Cu}(C_0(0,1] \otimes D_1))} \arrow[r] \arrow[d,
      "\textstyle \sigma_1"'] &
      {\text{Cu}(C_0(0,1] \otimes D_2))} \arrow[r] \arrow[d,
      "\textstyle \sigma_2"'] & \cdots \\
      \text{Cu}(B) \arrow[r
      ]                                              &
      \text{Cu}(B) \arrow[r
      ]                                              & \cdots 
\end{tikzcd}.
\]

The $*$-homomorphisms
$\pi_i : C_0(0,1] \otimes D_i \to B$ may be corrected by inner
automorphisms to obtain the following
one-sided approximate intertwining (in the sense of \cite[Section
2]{elliott}):
  \[
    \begin{tikzcd}
{C_0(0,1] \otimes D_1)} \arrow[r] \arrow[d, "\textstyle \pi_1"'] &
{C_0(0,1] \otimes D_2)} \arrow[r] \arrow[d, "\textstyle \pi_2"'] &
\cdots \arrow[r] & {C_0(0,1] \otimes D} \arrow[d, dotted, "\textstyle \pi"] \\
B \arrow[r]                                              & B
\arrow[r]                                              & \cdots
\arrow[r] & B
\end{tikzcd}
\]

By \cite[Remark 2.3]{elliott}, there exists a $*$-homomorphism
$\pi: C_0(0,1] \otimes D \to B$ as in the above diagram
and the map $\varphi : D \to B$ determined by the rule
$\varphi(d) := \pi (t \otimes d)$ is c.p.c.~ order zero, by
Corollary \ref{cone_oz}. Let us now verify
 that $\varphi$ satisfies the required tracial
 identity. By continuity, it
 suffices to check that the the tracial identity holds
at the finite stages. Let $d \in D_i$ and let $\tau \in TB$ be
given. Then 
\begin{align*}
  \tau \varphi^n(d)
  & \stackrel{(\ref{cone_eq})}{=} \tau \pi(t^n \otimes d)  =
    \tau \lim_{i \to \infty} \pi_i(t^n \otimes d) \\
  & =  \lim_{i \to \infty} (\mu \otimes (\iota_i^* \circ \Phi)(\tau))(t^n \otimes d) \\
  & = \mu(t^n) \Phi(\tau) (d).
\end{align*}
The second equality follows from the one-sided approximate
intertwining, and 
the third equality follows from continuity of $\tau$ and
by definition of $\sigma_i$.  (Note that
$\text{Cu}(\pi) = \lim_{i\to \infty} \sigma_i$,
but we do not actually use this.)
\end{proof}

\begin{lemma}
  There exists a sequence of fully supported
  Radon measures
  $\mu_k$ on $(0,1]$ with total mass one such that
  \begin{align}
    \lim_{k \to \infty} \mu_k(t) = 1.
    \label{leveled_moments}
  \end{align}
  Necessarily, this holds also for $t^n$, rather than just $t$,
  for each $n \in \mathbb N$.
  \label{measures}
\end{lemma}

\begin{proof}
  Pick measures which are increasingly weighted to the right.
\end{proof}


\begin{theorem}
  \label{thm}
  Any two unital, simple, $\mathcal Z$-stable, separable amenable
  C$^*$-algebras satisfying the UCT with
  isomorphic tracial simplices are coloured isomorphic.
\end{theorem}

\begin{proof} 
  We first consider the purely infinite case. This part of the
  theorem holds without assuming the UCT.
  Let $A$ and $B$ be Kirchberg algebras.
  By \cite[Corollary 2]{connes}, \cite[Theorem 3.1]{choi_effros},
  and \cite[Proposition 7]{wassermann}, separable amenable
  C$^*$-algebras are exact.
  Hence by Theorem \ref{o2-embedding}, there is a unital embedding
  of $A$ into
  $\mathcal O_2$. Since $B$ is purely
  infinite,  
  $K_0(B)$ consists of the Murray-von Neumann equivalence classes
  of non-zero properly
  infinite projections in $B$ (\cite[Theorem 1.4]{cuntz_k}).
  In particular, there is a
  non-zero properly infinite projection $p \in B$ whose $K_0$-class
  is equal to zero. It follows from \cite[Proposition 4.2.3]{rs} that
  there exists a unital embedding of $\mathcal O_2$ into $pBp$ and
  therefore a (possibly) non-unital embedding into $B$.
  Denote by $\varphi : A \to B$ the
  composition of these embeddings and construct $\psi : B \to A$
  symmetrically. Then since $\varphi$ and $\psi$ are injective and
  therefore isometric
  $*$-homomorphisms, the ultrapower maps 
  $\psi \varphi : A_\omega \to B_\omega$ and $\varphi \psi :
  B_\omega \to A_\omega$ induced by their compositions are also
  isometric and in particular injective.
  By two applications of Theorem \ref{kirchberg_colored},
  the pairs $\psi \varphi \iota_A$ and $\iota_A$;
  and $\varphi \psi \iota_B$ and $\iota_B$ are each coloured
  equivalent, where $\iota_A$ (resp.~ $\iota_B$) denote the canonical
  embedding of $A$ (resp.~ $B$) into its ultrapower.

  Now (as we may, in view of Theorem \ref{dichotomy}), suppose that $A$ and $B$ are two  finite C$^*$-algebras
  satisfying the hypotheses. By Theorem \ref{AF-embedding},
  there exists a unital embedding $\alpha_A$ of $A$ 
  into a separable, unital, simple AF algebra $D$ such 
  that the induced map $\alpha_A^*: T(D) \to T(A)$ is an
  isomorphism of tracial simplices.
  This, combined with an isomorphism $T(B)
  \to T(A)$ (assumed to exist), yields an isomorphism $\Phi : T(B)
  \to T(D)$. Pairing this with the faithful tracial states
  $\mu_k$ on 
  $C_0(0,1]$ of Lemma \ref{measures}, one obtains by Theorem
  \ref{induce_tensor_trace} a sequence of
  c.p.c.~order zero maps $\varphi_k : D \to B$ 
  satisfying the identity
  \begin{align}
    \tau \varphi_k^n = \mu_k(t^n) \Phi(\tau)
    \label{vf_tensor_trace}
  \end{align}
  for every $k,n \in \mathbb N$ and each $\tau
  \in T(B)$.
  Denote by $\varphi : A_\omega \to B_\omega$ the
  c.p.c.~order zero map induced between the ultrapowers by the maps $\varphi_k 
  \alpha_A : A \to B$. As above, by \cite[Corollary 3.4]{wz_cpc},
  $\varphi$ induces a
  mapping $\varphi^*$ of bounded traces on $B_\omega$ into  bounded traces on $A_\omega$.
  Then we have
  \begin{align}
    \begin{split} \label{vf_induce}
      \tau\varphi^n & \stackrel{(\ref{multiplicative_id})}{=}
      \lim_{k \to \omega} \tau_k (\varphi_k^n \alpha_A) \\
      &
      \stackrel{(\ref{vf_tensor_trace})}{=} \alpha_A^* (\lim_{k \to \omega} (\mu_k(t^n)
      \Phi(\tau_k)))
      \stackrel{(\ref{leveled_moments})}{=} \lim_{k \to \omega} (\alpha_A^* \Phi)(\tau_k)
    \end{split}
  \end{align}
  for each $n \in \mathbb N$ and limit trace $\tau = \lim_{k \to
    \omega} \tau_k \in
  T_\omega(B_\omega)$.

  Symmetrically, there is, again by Theorem \ref{AF-embedding}, a unital
  embedding $\alpha_B$ of $B$ into a unital, simple, separable AF algebra $E$
  such that the induced map $\alpha_B^* : T(E) \to T(B)$ is an
  isomorphism of tracial simplices. Using as above the isomorphism
  $\Psi := (\alpha_A^*  \Phi  \alpha_B^*) ^{-1} : T(A) \to
  T(E)$ and the faithful tracial states $\mu_k$,
  we obtain a c.p.c.~order zero map $\psi : B_\omega \to A_\omega$
  such that
  \begin{align}
    \tau\psi^n
    = \lim_{k \to \omega}(\alpha_B^* \Psi) (\tau_k)
    \label{psi_induce}
  \end{align}
  for each $n \in \mathbb N$ and limit trace
  $\tau = \lim_{k \to \omega} \tau_k\in T_\omega(A_\omega)$.
  So for all limit traces $\tau = \lim_{k \to \omega} \tau_k
\in T_\omega(A_\omega)$, we
  have 
  \begin{align} 
    \begin{split}
      \tau (\psi \varphi)^n
    & \stackrel{(\ref{multiplicative_id})}{=} \tau \psi^n
    \varphi^n 
     \stackrel{(\ref{psi_induce})}{=} (\varphi^n)^*
    (\lim_{k \to \omega}(\alpha _B^* \Psi)(\tau_k)) \\
    & \stackrel{(\ref{vf_induce})}{=}
    \alpha_A^* \Phi (\lim_{k \to \omega} (\alpha_B^* \Psi)(\tau_k))
    =    \lim_{k \to \omega} \alpha_A^* \Phi\alpha_B^*
    \Psi(\tau_k) 
    = \tau
    \end{split}
      \label{identitytrace} 
  \end{align}
  where the last equality follows from the definition of $\Psi$.

  Weak* density of $T_\omega(A)$ in 
  $T(A_\omega)$ (Theorem \ref{limit_traces}) extends 
  the identity (\ref{identitytrace}) to all tracial states on
  $A_\omega$. Since $A$ is simple, the canonical embedding of $A$
  into $A_\omega$ is totally full. 
  Hence by Theorem \ref{coloredclassification},
  $\psi \varphi \iota_A$ and $\iota_A$ are coloured equivalent.
  By a symmetric argument, $\varphi \psi \iota_B$ and $\iota_B$
  are coloured equivalent.
  It follows immediately from (\ref{vf_induce}) and
  (\ref{psi_induce}) that $\varphi$ and $\psi$ preserve constant
  limit traces.
  This shows $A$ is coloured isomorphic to
  $B$ (with $m$ and $n$ of (\ref{coloredid}) equal to two).
\end{proof}

\begin{corollary} Let 
  $A$ and $B$ be classifiable C$^*$-algebras (i.e., ones satisfying
  the hypotheses of Theorem \ref{thm}).
  Then the following statements are equivalent:
  \begin{enumerate}[(1)]
  \item There exist constant limit trace preserving
    c.p.c.~ order zero maps $\varphi :A_\omega \to
    B_\omega$ and $\psi : B_\omega \to A_\omega$ 
    such that $\varphi^n$ and $\psi^n$ induce mutually
    inverse isomorphisms of $T(A_\omega)$ and $T(B_\omega)$ for
    each $n \in \mathbb N$.
  \item There exist constant limit trace preserving
    c.p.c.~ order zero maps $\varphi :A_\omega \to
    B_\omega$ and $\psi : B_\omega \to A_\omega$ 
    such that $\varphi^n$ and $\psi^n$ induce mutually
    inverse isomorphisms of $T(A_\omega)$ and $T(B_\omega)$ for
    some $n \in \mathbb N$.
  \item $A$ and $B$ are coloured isomorphic (Definition \ref{coloured_isomorphic}).
  \item $A$ and $B$ are minimalist coloured isomorphic (Definition
   \ref{minimalist}).
  \item $T(A)$ is isomorphic to $T(B)$.
  \end{enumerate}
  \label{main_equivalence}
  Every isomorphism of $T(A)$ with $T(B)$ arises from a coloured
  isomorphism of $A$ and $B$.
\end{corollary}

\begin{proof}
  (1) $\Longrightarrow$ (2) is clear.
  (2) $\Longrightarrow$ (1) Let $\pi : A_\omega \to A^\omega$
  denote the canonical quotient map. By assumption,
  there exists an $n \in \mathbb N$ such that
  $\tau \psi^n
  \varphi^n = \tau$ for all $\tau \in T(A_\omega)$.
  By Remark \ref{uniform_tracial_ultrapower},
  this implies $\pi \psi^n
  \varphi^n = \pi \text{id}_{A_\omega}$. In particular, we will
  use that $\pi \psi^n \varphi^n(1_{A_\omega}) = \pi (1_{A_\omega})$
  below. Since 
  \[
    \pi \psi \varphi(1_{A_\omega})
    = \pi \psi \varphi(1_{A_\omega}^n)
    = (\pi \psi \varphi)^n(1_{A_\omega})
    = \pi \psi^n \varphi^n (1_{A_\omega}) = \pi (1_{A_\omega})
  \]
is a projection, it then follows from Corollary \ref{projection} that
$\pi \psi \varphi$ is a $*$-homomorphism. We have used Corollary
\ref{functional_calculus} and Corollary \ref{multiplicative} for
the second and third equalities. 
Therefore, $ \pi \psi^k \varphi^k = 
\pi (\psi\varphi)^k
= (\pi \psi \varphi)^k = \pi \psi \varphi$ for each $k \in \mathbb
N$. By Remark \ref{uniform_tracial_ultrapower} again,
this implies $\tau \psi^k \varphi^k = \tau$ for each $\tau \in
T(A_\omega)$ and each $k \in \mathbb N$. A symmetric argument
shows that $\varphi^k$ and $\psi^k$ induce mutually inverse
isomorphisms of $T(A_\omega)$ and $T(B_\omega)$.
  
  (1) $\Longrightarrow$ (3) follows (without the UCT
  assumption) from two
  applications of Theorem \ref{coloredclassification}: once with
  $\psi \varphi \iota_A$ and $\iota_A$ where $\iota_A$ is the
  constant sequence embedding of $A$ into $A_\omega$ and once more
  with $\varphi \psi \iota_B$ and $\iota_B$.
  (3) $\Longrightarrow$ (4) is immediate.
  (4) $\Longrightarrow$ (5) is a special case of Theorem
  \ref{preserve_traces}.
  (5) $\Longrightarrow$ (1) is given by (\ref{identitytrace}) in the
  proof of Theorem \ref{thm} (along with the analogous statement
  for $(\varphi \psi)^n$) and Corollary \ref{multiplicative}.

  The last statement follows from the
  proof of Theorem \ref{thm}, which constructs a coloured
  isomorphism using a given isomorphism of $T(A)$ with $T(B)$, and
  the proof of Theorem \ref{preserve_traces},
  which recovers the given
  isomorphism of $T(A)$ with $T(B)$ from the constructed coloured
  isomorphism. 
\end{proof}

  We aren't sure to what extent the c.p.c.~ order zero maps
  involved in a coloured isomorphism are unique.

\begin{corollary}
  \label{cor}
  Let $A$ be a classifiable C$^*$-algebra.
  If $T(A) = \varnothing$, then $A$ is
  coloured isomorphic 
  to $\mathcal O_2$. If $T(A) \neq \varnothing$, then $A$ is
  coloured 
  isomorphic to a unital simple AF algebra.
\end{corollary}

\begin{proof}
  If $T(A) = \varnothing$,
  then by Theorem \ref{thm} $A$ is coloured isomorphic to
  $\mathcal O_2$.
  If $T(A) \neq \varnothing$, 
  by \cite[Theorem 3.10]{blackadar}, there
  exists a unital AF algebra $B$ whose tracial simplex 
  isomorphic to $T(A)$. The conclusion now
  follows from Theorem \ref{thm}. 
\end{proof} 

The existence step in establishing finite nuclear dimension from
$\mathcal Z$-stability (\cite[Lemma 7.4]{bbstww} and \cite[Lemma
5.2]{cetww}) is the construction of a sequence of
c.p.c.~ maps $\phi_i : A \to A$,
where $A$ is a finite, unital, simple, $\mathcal Z$-stable,
separable amenable C$^*$-algebra, 
which factorize through
finite-dimensional C$^*$-algebras $F_i$ as
\[
\begin{tikzcd}
  A \arrow[rr, "\textstyle \phi_i"] \arrow[rd, "\textstyle \theta_i"'] &                        & A \\
  & F_i \arrow[ru, "\textstyle \eta_i"'] &  
\end{tikzcd}
\]
with $\theta_i$ c.p.c.~ and $\eta_i$ c.p.c.~ order zero, such that
the induced maps $(\theta_i)_{i=1}^\infty : A \to \prod_\omega
F_i$ and $\Phi = (\phi_i)_{i=1}^\infty : A \to A_\omega$ are order
zero and
\[
  \tau \Phi(a) = \tau(a)
\]
for each $a \in A$ and each $\tau \in T(A_\omega)$.

More precisely, the conclusion in \cite[Lemma 5.2]{cetww} is that
$\Phi$ agrees 
with the canonical embedding $\iota_A$ of $A$ into $A_\omega$ in
the uniform 
tracial ultrapower $A^\omega$ 
while the conclusion in \cite[Lemma 7.4]{bbstww} is that $\tau
\Phi = \tau \iota_A$ for each $\tau \in T(A_\omega)$. These
conclusions are equivalent by Remark
\ref{uniform_tracial_ultrapower} and Theorem \ref{limit_traces}.

In the presence of the UCT,
a one-sided formulation of Theorem \ref{thm} (Corollary
\ref{existence_step} below) can be viewed as a
generalization of the existence step mentioned earlier
since the conclusion follows from the special case that $B$ is an
AF algebra.
Since the nuclear dimension of $B$ is zero
(\cite[Remark 2.2 (iii)]{wz_nd}), there exist finite-dimensional
C$^*$-algebras $F_i$ and c.p.c.~ maps
$\rho_i : B \to F_i$ and c.p.c.~ order zero maps $\sigma_i :
F_i \to B$ such that the triangle below commutes approximately.
\[
  \begin{tikzcd}
    A \arrow[rr, "\textstyle \phi_i"] \arrow[d, "\textstyle \varphi_i"']        &                                        & A                                 \\
    B \arrow[rr, "\textstyle \text{id}_B"] \arrow[rd, "\textstyle \rho_i"'] &                                        & B \arrow[u, "\textstyle \psi_i"'] \\
    & F_i \arrow[ru, "\textstyle \sigma_i "'] &                                  
  \end{tikzcd}
\]
The maps $\theta_i: A \to F_i$ on the left side of the diagram
and the maps $\eta_i : F_i \to A$ on the right side of the diagram
satisfy the conclusion of the aforementioned existence step.

\begin{corollary}\label{existence_step}
  Let $A$ and $B$ be finite classifiable C$^*$-algebras
  with isomorphic tracial simplices.
  Then there exist sequences of c.p.c.~ order
  zero maps $\varphi_i : A \to B$ and $\psi_i : B \to A$, $i \in
  \mathbb N$, such that the induced c.p.c.~ order zero maps
  $\varphi : A_\omega \to B_\omega$ and $\psi : B_\omega \to A_\omega$ induce mutually
  inverse isomorphisms of $T(A_\omega)$ and $T(B_\omega)$. In
  particular,
\[
  \tau \psi \varphi(a) = \tau(a)
\]
for each $a \in A$ and each $\tau \in T(A_\omega)$.
\end{corollary}

\begin{proof}
  This is a special case of Corollary \ref{main_equivalence}.
\end{proof}

Corollary \ref{existence_step} gives rise to the
completely positive approximation property below.
The nuclear dimension calculation in Corollary \ref{nuclear_dim}
is not new -- in fact, it 
holds without the UCT 
assumption and without unitalness (\cite{cetww}, \cite{ce}), and
our proof of it using the UCT still relies on the
main technical results of \cite{bbstww} and \cite{cetww} that were
used in establishing finite nuclear dimension from $\mathcal
Z$-stability in the context of the Toms-Winter conjecture
(\cite{cetww}, \cite{ce}).

\begin{corollary} \label{nuclear_dim}
  Let $A$ and $B$ be finite classifiable
  C$^*$-algebras with isomorphic tracial simplices.
  Then there exist a sequence
  $\varphi_i : A \to B$ of c.p.c.~ order zero maps and a sequence $\xi_i : B \to
  A$ such that $\xi_i$ is a sum of two c.p.c.~ order zero
  maps from $B$ to $A$, $i \in \mathbb N$, and such that the
  following diagram approximately commutes:
  \[
    \begin{tikzcd}
      A \arrow[rr, "\textstyle \mathrm{id}_A"] \arrow[rd, "\textstyle \varphi_i"'] &                                   & A \\
      & B \arrow[ru, "\textstyle \xi_i"'] &  
    \end{tikzcd}.
  \]
  Since a possible choice for $B$ is an AF algebra, it follows
  that the nuclear dimension of $A$ is at most one. If $A$ is not
  an AF algebra, then the nuclear dimension of $A$ is exactly one.
\end{corollary}

  \begin{proof}
    Let notation for $\varphi_i$, $\psi_i$, $\varphi$, and $\psi$
    be as in Corollary
    \ref{existence_step}, and let $h$ be a positive
    element with full spectrum in $\mathcal Z$. Then since
    $\iota_A$ is totally full and since $\iota_A$ and $\psi
    \varphi \iota_A$ agree on traces, $\iota_A \otimes h$ and
    $\psi \varphi \iota_A \otimes h$ are approximately unitarily
    equivalent, by \cite[Lemma 4.8]{cetww}.
    Since $1_{\mathcal Z} - h$ is also a positive
    element of $\mathcal Z$ with full spectrum, $\iota_A \otimes
    (1_{\mathcal Z} - h)$ and $\psi \varphi \iota \otimes
    (1_{\mathcal Z} - h)$ are also approximately unitarily
    equivalent.
    Therefore, there exist unitaries $u_1$ and
    $u_2 \in (A \otimes \mathcal Z)_\omega$ such that
    $a \otimes h = {u_1} (\psi \varphi(a) \otimes h) u_1^*$ and
    $a \otimes (1_{\mathcal Z} - h) = 
    {u_2} (\psi\varphi(a) \otimes (1_{\mathcal Z} -
    h)) u_2^*$.
    Since $\mathcal Z$ is strongly
    self-absorbing (\cite[Theorem 8.7]{jiang-su}) 
    and therefore $K_1$-injective (\cite[Theorem 6.7]{rordam},
    \cite[Theorem 10.12]{rieffel}),
    there is a
    $*$-homomorphism $\theta : A \otimes \mathcal Z \to
    A_\omega$ such that $\theta(a \otimes 1_{\mathcal Z}) = a$
    for each $a \in A$, by \cite[Theorem 2.3]{toms-winter}.
    Let $h_1 := h$ and $h_2 := 1_{\mathcal Z}-h$.
    Taking representative sequences of $*$-homomorphisms
    $\theta_i : A \otimes 
    \mathcal Z \to A$ (which exist by the proof of \cite[Theorem
    2.3]{toms-winter}) and of unitaries
    $(u_1^{(i)})_{i=1}^\infty$ and $(u_2^{(i)})_{i=1}^\infty$ in $A
    \otimes \mathcal Z$
    corresponding to $u_1$ and $u_2$, we have
    the following approximately commutative diagram.
\[    
      \begin{tikzcd}
        A \arrow[r, "\textstyle \mathrm{id}_A"] \arrow[rd, "\textstyle \varphi_i"'] & A                                 & A \otimes \mathcal{Z} \arrow[l, "\textstyle \theta_i"']                          \\
        & B \arrow[r, "\textstyle \psi_i"'] & A \arrow[u, "\textstyle \sum_{k=1}^2 u_k^{(i)} (\cdot \otimes h_k) u_k^{(i)*}"']
      \end{tikzcd}
    \]
    The maps $\varphi_i : A \to B$ and the maps $\xi_i : B \to A$
    which factor though $A$ 
    and $A \otimes \mathcal Z$ give the desired completely positive
    approximation property.

    Let us now specialize to the case that $B$ is an AF algebra
    with $T(B) \cong T(A)$, which exists by \cite[Theorem
    3.10]{blackadar}.
    Since the nuclear dimension of $B$
    is zero (\cite[Remark 2.2 (iii)]{wz_nd}),
    there exist c.p.c.~ maps
    $\rho_i : B \to F_i$ into finite-dimensional C$^*$-algebras
    $F_i$ and c.p.c.~ order zero maps $\sigma_i :
    F_i \to B$ which make the following diagram approximately
    commute
\[
      \begin{tikzcd}
A \arrow[d, "\textstyle \varphi_i"'] \arrow[rr, "\textstyle \text{id}_A"] &                                        & A                                \\
B \arrow[rd, "\textstyle \rho_i"'] \arrow[rr, "\textstyle \text{id}_B"]   &                                        & B \arrow[u, "\textstyle \xi_i"'] \\
                                                                    & F_i \arrow[ru, "\textstyle \sigma_i"'] &                                 
\end{tikzcd}.
\]
Since the maps $\rho_i \varphi_i : A \to F_i$ are 
and the maps $\xi_i \sigma_i : F_i \to A$ are c.p.~ and since
$\xi_i \sigma_i$ is a sum of two c.p.c.~ order zero maps,
 $A$ has nuclear dimension at most one (\cite[Section
1.1]{cetww}). If $A$ is not an AF algebra, then the nuclear
  dimension of $A$ is exactly one (\cite[Remark 2.2
  (iii)]{wz_nd}).
\end{proof}

\begin{question}
  All that is used about the (multiplicative)
  AF embeddings in Theorem \ref{thm} is
  that they are c.p.c.~ order zero maps which induce
  isomorphisms of simplices. In fact, it would be 
  enough for each AF embedding to be replaced by
  a sequence of order zero maps which induce
  an isomorphism of tracial
  simplices in the sense described in Section
  \ref{traces}. If this could be done without the
  UCT, then one would have Theorem \ref{thm} for classifiable
  C$^*$-algebras not necessarily satisfying the UCT, as the UCT
  assumption is only used to produce AF embeddings.
  If $A$ and $B$ are unital, simple,
  exact, $\mathcal Z$-stable, separable C$^*$-algebras
  with stable rank one,
  is there an order zero map $\varphi : A \to B$ realizing
  prescribed tracial data as in Theorem \ref{induce_tensor_trace}? 
\end{question}

\bibliography{mybib.bib} 

\vfill

\noindent  \textsc{Department of Mathematics,
  University of Toronto,  \\
  40 St. George Street, Toronto, ON, Canada
 \hspace{0 pt} M5S 2E4
}\par\nopagebreak
\noindent   \textit{E-mail addresses}: \texttt{elliott@math.toronto.edu}, \texttt{jim@math.toronto.edu}

\end{document}